\documentclass{amsart}

\usepackage{amsfonts,amsthm,amssymb,latexsym,}
\usepackage{mathrsfs}
\usepackage[all]{xy}
\usepackage[dvips]{graphicx}
\usepackage[latin1]{inputenc}

\newtheorem{thm}{Theorem}[section]
\newtheorem{prop}{Proposition}[section]
\newtheorem{lem}{Lemma}[section]
\newtheorem{cor}{Corollary}[section]
\newtheorem{rem}{Remark}[section]
\newtheorem{example}{Example}[section]

\newcommand{\R}{I\!\! R}
\newcommand{\N}{I\!\! N}

\newcommand{\Ttheta}{T_{\theta} TM}
\newcommand{\J}{\mbox{J}}

\newcommand{\dt}{\frac{d}{dt}}
\newcommand{\ds}{\frac{d}{ds}}
\newcommand{\Dt}{\frac{D}{dt}}
\newcommand{\KM}{K_{\mbox{mag}}}
\newcommand{\Dtt}{\frac{D^2}{d^2t}}

\newcommand{\Z}{Z\!\!\! Z}

\newcommand{\ps}{\frac{\partial}{\partial s}}

\newcommand{\m}{\mathcal}

 \def \co{{\mathbb C}}
 \def \di{{\mathbb D}}

\begin{document}

\title{ Positive topological entropy  for magnetic flows on surfaces. }

\author{Jos\'e Ant\^onio Gon\c{c}alves Miranda}

\address{ Universidade Federal of Minas Gerais, Dep. of Matem\'atica,
ICEx \\
Av. Ant\^onio Carlos, 6627/ C. P.702, 30123-970, Belo Horizonte,
MG, Brasil.}
 \email{jan@mat.ufmg.br}

% 37B40 Topological entropy
% 37D30 Partially hyperbolic systems and dominated splittings
% 37J50 Action-minimizing orbits and measures
% 37J40 Perturbations, normal forms, small divisors, KAM theory, Arnold diffusion
% 37Jxx Finite-dimensional hamiltonian, lagrangian, contact, and nonholonomic systems
% 70K45 Normal forms
\begin{abstract}
We study the topological entropy of the magnetic flow on a closed
riemannian surface. We prove that if the magnetic flow has a
non-hyperbolic closed orbit in some energy set $ T^cM= E^{-1}(c)
$, then there exists an exact $ C^\infty$-perturbation of the
2-form $ \Omega $   such that the new magnetic flow has positive
topological entropy in $ T^cM $. We also  prove that if the
magnetic flow has an infinite number of closed orbits in $ T^cM $,
then there exists an exact $ C^1 $-perturbation of $ \Omega $ with
positive topological entropy in $ T^cM $. The proof of the last
result is based on an analog of Franks' lemma for magnetic flows
on surfaces, that is proven in this work, and  Mañé's techniques
on dominated splitting. As a consequence of those results, an
exact magnetic flow on $ S^2$ in high energy levels admits a $ C^1
$-perturbation with positive topological entropy. In the
appendices we show that an  exact magnetic flow on the torus  in
high energy levels admits a $ C^\infty $-perturbation with
positive topological entropy.
\end{abstract}

\maketitle

\section{Introduction and statements.}
Let $ M $ be a closed and  oriented surface with a smooth
riemannian metric $ g $ and $ \pi: TM \rightarrow M $ its tangent
bundle. Let $ \omega_0 $ be the symplectic structure on  $ TM $
obtained  by  pulling back  the canonical symplectic structure of
the cotangent bundle $ T^*M $ using  the riemannian metric $ g$.
For any  $ x \in M $, let $ i: T_xM\rightarrow T_xM $ be the
linear map such that $ \{{v}, i \cdot v\} $ is a positive oriented
orthogonal basis for $ T_xM$. Consider the 2-form $ \Omega_0 $ in
$ M $ defined by:
$$ \Omega_0(x)( u,v) = g_x(\ i \cdot u\ ,\ v\ )\ \ \mbox{(  Area form )}. $$
We denote by $ \Omega^2(M) $ the set of all smooth 2-forms on $
M$. Since $ M $ is a surface, if $ \Omega \in \Omega^2(M) $, there
exist a smooth function $ f:M \rightarrow \R $ such that $ \Omega=
f\cdot \Omega_0$. Therefore, we can consider the $ C^k$-topology
in $ \Omega^2(M) $ as being the $ C^k$-topology in the space of
smooth functions on $ M $, which we denote by $ C^k (M)$.

\smallskip

 Given $ \Omega = f \cdot \Omega_0\in \Omega^2(M)$, let $
\omega(\Omega) $ be the  symplectic structure in $ TM $ defined by
 $$\omega(\Omega) =
\omega_0 + \pi^* \Omega= \omega_0 + ( f\circ \pi)\cdot  \pi^*
\Omega_0
$$
 which is called the {\it twisted  symplectic structure}.
Let $ E: TM \rightarrow \R $ be the hamiltonian given by
$$E (x,v) = \frac{1}{2}  g_x(v,v)\ \ \mbox{
(Kinetic Energy)}.$$The {\it magnetic field} associated to $
\Omega  $ is the hamiltonian field $ X(\Omega)=X^f $ of the
hamiltonian $ E $ with respect to $ \omega(\Omega)$. The {\it
magnetic flow} associated to $ \Omega $ is the hamiltonian flow $
\phi_t(\Omega)=\phi_t^f: TM\rightarrow TM $ induced by the vector
field $ X(\Omega)$. This flow models the motion of a unit mass
particle under the effect of the {\it Lorentz force} $ Y= f\cdot i
$. In other words,  a curve $t \mapsto  ( \gamma(t),
\dot\gamma(t)) \subset TM $ is an orbit of $ \phi_t$  if and only
if $ \gamma: \R \rightarrow M $ satisfies:
\begin{equation}
\label{eq1}
 \frac{D}{dt} \dot\gamma = Y_\gamma( \dot\gamma) = f(\gamma) \ i
 \cdot \dot \gamma .
\end{equation}
Observe that if $ \Omega \equiv 0 $ (i.e. $ f\equiv 0 $), the
above equation coincides with the geodesic equation on $ M $ for
the riemannian metric $g$.   A curve which satisfies  equation
(\ref{eq1}) is called an {\it $ \Omega$-magnetic geodesic}.

For $ c>0$,  let $ T^cM $ be the bundle defined by $ T^cM = E^{-1}
(c) $. Note that $ T^cM $ is a compact invariant  submanifold of $
TM $ and  that the restriction of $ \phi^\Omega_t $ to $ T^cM $
has no fixed points.

 When   $ \Omega $ is  an exact 2-form (i.e.
$ \Omega= d\eta $), we can define  the convex and superlinear
lagrangian $L_\eta :TM \rightarrow \R $ as
   $$ L_\eta(x,v) = \frac{1}{2}g_x(v,v) - \eta_x(v).$$
Computing the  Euler-Lagrange equation of $ L $, we obtain that
the  extremal curves coincide with the $d\eta$-magnetic geodesics.
Then the magnetic flow associated to $ d\eta $  can be studied as
a lagrangian flow. In this case the flow is called {\it exact
magnetic flow}.

\smallskip

Magnetic flows  have attracted considerable attention in recent
years. This class of dynamical systems was first considered by
V.I. Arnold in \cite{Arnold:P61} and by D. V. Anosov and Y. G.
Sinai in \cite{Anosov-Sinai:P-67}. remarkable properties of
magnetic flows were obtained by many authors; we refer to
\cite{Ginzb:p96, P-P:c96, Mario:p99, Burns-Part:p2002,
C-M-p:p2004}. In this work we are interested to study  the
behavior of the functional $ \Omega\mapsto h_{top} ( \Omega, c )$
for a prescribed  energy level $ c>0$, where $ h_{top} ( \Omega, c
)= h_{top} ( g,\Omega, c ) $ denotes the topological entropy of $
\left.\phi^\Omega_t\right|_{T^cM} $.

The {\it topological entropy} is a dynamical invariant that,
roughly speaking, measures its  orbit structure complexity. Its
precise definition can be found in \cite{Bow:b75}. The relevant
question about the topological entropy is whether it is positive
or vanishes.   Standard arguments in dynamical systems show that
if a flow contains a transversal homoclinic orbit then it has
positive topological entropy (in fact it contains a nontrivial
hyperbolic set). Conversely, if $ h_{top} ( \Omega, c )> 0 $, the
result of A. Katok \cite{Kat:p80} implies  that   $ \phi^\Omega_t
$ presents a transversal homoclinic orbit in $ T^cM$. In
particular, it has infinitely many closed orbits.

 The main results of this paper are:
\begin{thm}\label{T_maim-e}
 Let $ M $ be a closed oriented  surface with a smooth metric $ g$. Let  $ \Omega
$ be a  2-form on M and  $ c > 0$.  Suppose that the magnetic flow
$ \phi_t^\Omega $ has a non-hyperbolic closed orbit in    $ T^cM$.
Then there is an exact 2-form $ d\eta $ of arbitrarily small norm,
in the $C^r$-topology  (with $ 4 \leq r \leq \infty $), such that
$ h_{top}(\Omega+d\eta,c)> 0 $.
\end{thm}
\begin{thm}\label{T_maim-h}
 Let $ M $ be a closed oriented  surface with a smooth metric $ g$. Let  $ \Omega
$ be a  2-form on M and  $ c > 0$. Suppose that the  magnetic flow
$ \phi_t^\Omega $ has an infinite number of  closed orbits in    $
T^cM $. Then there is an exact 2-form $ d\eta $ of arbitrarily
small norm,  in the $C^1$-topology,  such that $
h_{top}(\Omega+d\eta,c)> 0 $.
\end{thm}

Two important tools to prove these theorems are a version of the
conservative Kupka-Smale theorem for  magnetic flows on surfaces,
and a family of generic properties for the $k$-jet of the Poincaré
map of all  closed orbits, both proven in
\cite{jan-generic:p2006}. We will give the precise statements in
section 2.  In order to prove theorem \ref{T_maim-e}, we will show
that if there exists a non-hyperbolic closed orbit in $ T^cM $,
then   we can approximate   $ \Omega $ by another  2-form in the
same cohomology class of $ \Omega$, such that, the corresponding
magnetic flows has an elliptic closed orbit in $ T^cM $ for which
the Poincaré map is  a generic exact twist map in a neighborhood
of the elliptic fixed point. Then, a result of Le Calvez
\cite{Calv:p91} implies that the Poincaré map has a transversal
homoclinic point. Therefore, the magnetic flow can be approximated
by another one with positive topological entropy. We will give the
details of these arguments in section 3. Hence, in order to prove
theorem \ref{T_maim-h}, we can assume that all closed magnetic
geodesics in $ T^cM $ are hyperbolic. Using Mañé's technique of
dominated splitting \cite{Man:p82} and an analog of the Franks'
lemma (theorem \ref{franks-L}) for magnetic flows on surface, %that
%we will prove in  section \ref{Fraks-sec1},
we will obtain a nontrivial hyperbolic set of $ \phi_t^\Omega $ in
$ T^cM$. Since Mañé's technique and  Franks' lemma only work in
the $ C^1$-topology,  we  only obtain this approximation in the $
C^1$-topology. The details and statements are given in section
\ref{th-h}.

 \smallskip

Let us now describe an application of theorems \ref{T_maim-e} and
\ref{T_maim-h} to exact magnetic flows on $ S^2$ in high energy
levels. We recall the  definition of the  strict Mañé's critical
value for  convex and superlinear lagrangians (cf.
 \cite{Man:c96,C-D-I:p97} and \cite{P-P:p97}).
Let  $ L : TM \rightarrow \R $ be a convex and superlinear
lagrangian. The {\it strict Mañé's critical values} of $ L $  is
the real number $ c_0(L) $ such that
$$c_0(L)= \inf\{k\in \R; \ \int_0^T \ (L(\gamma(t),\dot\gamma(t))+k) \
dt \geq 0 \ \mbox{ for any absolutely continuous closed curve }
\gamma $$ $$ \mbox{
 homologous to zero defined in any closed interval [0,T]}\}.$$
It is well known that for  an arbitrary surface $ M $, if $ \Omega
=
 d\eta $ and $ c > c_0(L_\eta ) $, then
the  restriction of the  exact magnetic flow  in the   energy set
$ T^cM $ is a  reparametrization  of a  geodesic flow  in the unit
tangent bundle for an  appropriate   Finsler metric on $ M$  (cf.
\cite{C-I-P-P:p98}). Recall that a  {\it  Finsler metric} is a
function $ F:TM\rightarrow \R $ such that: $ F $ is differentiable
away from the zero section, its second derivative in the direction
of the fibres is positive defined and $ F(x,\lambda\ v ) = \lambda
F( x,v ) $ for all $ \lambda
> 0 $ and $ (x,v)\in TM $.  If $
g $ is a  riemannian metric  on $ M$, then $ F(x,v) =
g_x(v,v)^{1/2} $ is a trivial  example of a  Finsler metric. We
say that  a  Finsler metric is {\it bumpy } if  all closed
geodesics are {\it non-degenerate}, this is, if the linearized
Poincaré map of every closed geodesic does not admit a root of
unity as an eigenvalue.

  Many results for geodesic flows
of a riemannian metric  remain valid  for Finsler metrics, but, in
contrast with the  riemannian case, there exist examples of bumpy
Finsler metrics on $ S^2$ with only  two closed geodesics. These
examples were given by  Katok in \cite{Kat:p73} and were studied
geometrically by Ziller in \cite{Ziller:p83}.

The following  theorem is a particular case of a  result proved by
Radamacher in \cite[theorem 3.1(b)]{Rad:p89} for bumpy geodesic
flows on compact simply-connected manifold satisfying a
topological condition over its rational cohomology algebra  $
H^*(M,\mathbb{Q})$. This condition holds  for  $ S^2$ and this
result remains valid for bumpy Finsler metrics (cf. \cite[ pg.
81]{Rad:p89}). See also the theorem proved by N. Hingston in
\cite[section 6.2]{Hingston:p84}.

\begin{thm} \label{T-Rad} Let $ F:TS^2\rightarrow \R
$ be a bumpy  Finsler metric on $ S^2 $. Suppose  that there  are
only finitely many  closed geodesics for  $ F $ in $ S^2 $. Then
there is least one  non-hyperbolic closed geodesic.
\end{thm}

Combining  theorem \ref{T-Rad} with   theorems \ref{T_maim-e} and
\ref{T_maim-h}, we obtain the following proposition, which is a
version of a result for geodesic flows on $ S^2 $ proved by G.
Contreras and G. Paternain in \cite{C-P:p2002}, for the class of
flow studied here.
\begin{prop}\label{c-pert-c1-S}
Let  $ \Omega= d\eta $ be an exact 2-form in $ (S^2,g) $.  Then,
for any $ c >c_0(L_\eta)$
 there is an  exact 2-form  $ d\overline \eta $ on $S^2 $, of
norm arbitrarily small in the $ C^1$-topology, such that $
h_{top}(d(\eta+\overline \eta),c) > 0 $.
\end{prop}

A 2-form $ \Omega $ on $ M$ is said to be {\it weakly exact} if
its lift to the universal covering $ \tilde M$  of $ M  $ is
exact. Of course an exact form is weakly exact. If $ \Omega $ is
weakly exact then the lift of the magnetic flow to $ \tilde M $ is
an exact magnetic flow and we can define the critical value $
c(g,\Omega)=c(\Omega)$ as the strict critical value of the
lagrangian on $ T \tilde M $ corresponding to the lifted flow on $
T \tilde M $, that  can be infinite. In fact $ c(\Omega)< \infty $
if and only if the lift of $ \Omega $  has a bounded primitive
(cf. \cite{Burns-Part:p2002}). For  surfaces $ M $ of genus $ \geq
2$, each  $ \Omega\in \Omega^2(M) $ is weakly exact and we always
have $ c(\Omega)<\infty$. In this case, K. Burns and G.Paternain
proved that $ h_{top}(\Omega,c)
> 0 $ for all  $ \Omega \in \Omega^2(M)$ and for all $ c >
c(\Omega) $ \cite[proposition 5.4]{Burns-Part:p2002}.

In   Appendix \ref{A-torus}, using   results of \cite{Math:p91},
\cite{C-P:2002-3} and \cite{Mass:p97}, we  prove a result
analogous to proposition \ref{c-pert-c1-S} for the two-dimensional
torus by performing  perturbations in the $ C^\infty $-topology.

\section{preliminaries: generic properties}

In this section we will give  the statements of some results
proved in \cite{jan-generic:p2006} that we shall use in the proof
of the main results of this work.

 We say that a property $ P $ is {\it $ C^r$-generic for
 magnetic flows} if, for any  $ c > 0 $,
   there exists a   subset $ \m O(c) \subset
\Omega^2(M) $, such that:
\begin{itemize}
\item[\rm (a)] The subset
$\mathcal O_h(c) := \mathcal O(c) \cap \{ \Omega \in \Omega^2(M) ;
[\Omega]= h \} $ is $ C^r$-residual in $ \{ \Omega \in \Omega^2(M)
; [\Omega]= h \} $, for all $ h \in H^2(M,\R)$.
\item[\rm (b)]The flow $ \left.\phi_t^\Omega \right|_{T^cM} $ has the property $
P $, for all $ \Omega \in \m O(c)$.
\end{itemize}

The following theorem is a conservative version of the Kupka-Smale
theorem for magnetic flows on surfaces.

\begin{thm}
\label{T-KS}\cite[theorem 1.2]{jan-generic:p2006} Let $ M $ be a
closed and oriented surface with a smooth metric $ g $. The
following property:
$$ P_{K-S}:\left\{\begin{array}{clcr}
 (i)&\mbox{all closed orbits are
hyperbolic or  elliptic}\\
 (ii)&\mbox{all  heteroclinic points are   transversal}
 \end{array}\right.$$
 is $ C^r$-generic for  magnetic flows on  surfaces,  with  $ 1\leq r \leq \infty$.
 \end{thm}
Let us  recall some facts about  the jet space for symplectic maps
in $ (\R^{2n},\omega_0)$. Let $ Diff_{\omega_0}(\R^{2n},0)$ be the
space of smooth symplectic diffeomorphisms $ f :(\R^{2n},\omega_0)
\rightarrow (\R^{2n},\omega_0) $ that fix the origin. Given $ k
\in \N $, consider the equivalence relation
 $ \sim_k $ in $ Diff_{\omega_0}(\R^{2n},0)$, defined as:
$$ f \sim_k g \Leftrightarrow  \mbox{the Taylor polynomials of degree } k
\mbox{ at zero are equal}.$$ We define the {\it k-jet } of $ f \in
Diff_{\omega_0}(\R^{2n},0)$, which we will  denote by $
j^k(f)=j^k(f)(0) $, as the equivalence class with respect to the
relation $ \sim_k$. The {\it space of symplectic k-jets} $
J^k_s(n) $ is the set of all equivalence class with respect to the
relation $ \sim_k $ of elements of $ Diff_{\omega_0}(\R^{2n},0)$.
When  $ k = 1 $, we can identify  $ J^1_s(n) $ with  $ Sp(n)$. We
say that a subset
 $ Q \subset  J^k_s(n)$ is {\it invariant} when
 $ \sigma \cdot Q \cdot \sigma^{-1} =
Q $, for all $ \sigma \in J^k_s(n) $.

Let $ \theta_t = \phi_t^\Omega (\theta ) $ be a periodic orbit
with period $ T> 0 $ in $T^cM$ and  $ \Sigma \subset  T^cM$ be a
local transversal section in the energy level $ T^cM  $ over the
point $ \theta $. Then, the twisted symplectic form $
\omega(\Omega) $ induces a symplectic form on $ \Sigma $ and the
Poincaré map $ P(\theta,\Sigma,\Omega):\Sigma \rightarrow \Sigma $
preserves this form. Therefore, using Darboux coordinates, we can
assume that $ j^k( P(\theta,\Sigma,\Omega))\in J^k_s(1)$. The fact
that $ j^k( P(\theta,\Sigma,\Omega))$ belongs to an invariant
subset $ Q $ is independent of the chosen section $ \Sigma\subset
T^cM $ and of the chosen coordinates  of $ \Sigma$.
\begin{thm}
\label{T-jet}\cite[theorem 1.3]{jan-generic:p2006}
%Given a open
%and dense invariant subset $ Q \subset J^k_s(1) $,
%  the property $$P_{Q} :
%\mbox{ the k-jet of the  Poincar\'e map belong  to  }  Q $$
% is
%a $ C^r$-generic property for  Magnetic Flows on surfaces,
% with $ k < r \leq \infty$.
Let $ \theta_t= \phi_t^\Omega(\theta) $ be a closed orbit. Let $ Q
\subset J^k_s(1) $ be an open and invariant, such that $j^k(
P(\theta,\Sigma,\Omega)) \in \overline Q$. Then there exists an
exact 2-form  $ d \eta \in \Omega^2(M) $, arbitrarily  $
C^r$-close to zero, with $ r
> k$, such that
\begin{itemize}
\item[\rm(i)]
$ \theta_t$ is a closed orbit of $ \phi_t^{\Omega+d\eta} $ and
\item[\rm(ii)]
 $ j^k( P(\theta,\Sigma,\Omega+d\eta) ) \ \in \ Q$.
\end{itemize}
\end{thm}

\section{magnetic flows with a non-hyperbolic closed orbit}

%In this section we will prove the theorem \ref{T_maim-e}.

Let us  recall the Birkhoff's Normal  Form (for a proof  see
\cite[pg. 222]{Sieg-Mos:b95}).

\begin{thm}\label{FNB}
Let f  be a $ C^4 $ diffeomorphism defined in a neighborhood of $
0 \in \R^2 $ such that  $ f $ preserves the area  form  $ dx\wedge
dy $ and $f(0)=0$.  Suppose that the eigenvalues of  $ d_0f $
satisfy: $ | \lambda | = 1 $ and $\lambda^n \not= 1$, for all $ n
\in \{1,...,4\}$. Then there exists  a $ C^4 $ diffeomorphism $ h
$, defined in a neighborhood of $ 0 $ such that:  $ h(0) =0 $, $ h
$ preserves the form $ dx\wedge dy $ and
%$$h^{-1} \circ f \circ h
%\left(\begin{array}{cc} x\\y
%\end{array}\right) =\left( \begin{array}{crcl}
%\cos( \alpha + \beta( x^2+ y^2)) & - \sin ( \alpha + \beta( x^2+
%y^2))\\  \sin ( \alpha + \beta( x^2+ y^2))& \cos  ( \alpha +
%\beta( x^2+ y^2))\end{array}\right)\left(\begin{array}{cc} x\\y
%\end{array}\right) + \mathcal O ( ( x^2+ y^2)^2)
%$$ where $ \lambda = e^{\pm 2\pi i \alpha}$  is eigenvalue of $ d_0f $.
%
%the map $ h^{-1} \circ f \circh $ is written as:
%we have:
$$ h^{-1} \circ f \circ h ( r, \theta ) = ( \ r \ , \ \theta +
\alpha+ \beta r^2 \ ) + \mathcal O( r^4)$$ in polar coordinates
$(r,\theta)$. Moreover, the property of $ \beta\not= 0 $  depends
only on   $ f$.
\end{thm}

We say that a  homeomorphism $ f : [ a,b ] \times S^1 \rightarrow
[ a,b] \times S^1 $ is a  { \it twist map} if for all  $ \theta
\in S^1 $ the function $ [a,b] \mapsto \pi_2\circ f(\cdot ,
\theta) \in S^1$ is strictly  monotonic . Observe that if the
coefficient $ \beta= \beta(f) $ in the normal form  is not equal
to zero, then
% for $ |r|\leq \epsilon$ with $\epsilon$ small
%enough,
$ f $ is conjugated to a  twist map in $ [0,\epsilon]\times S^1$,
for $\epsilon$ small enough.

 We shall use the following result:
\begin{prop}[Le Calvez~{\cite[Remarques pg.34]{Calv:p91}}]
\label{p-calves}
 Let $f$ be a diffeomorphism of the annulus $\R\times S^1$ that
  is an area preserving twist map  and  is such that the form
 $f^*(R\,d\theta)-R\,d\theta$ is exact. Suppose  $ f $ has the following  properties:
  \begin{itemize}
 \item[{(i)}] for every periodic point  $x$  of period  $q \geq 1 $,
the real number $1$ is not an  eigenvalue of $d_xf^q$,
 \item[{ (ii)}] the stable and unstable manifolds of any couple of  hyperbolic
           periodic points of $f$ intersect transversally
           (i.e. whenever they meet, they meet transversally).
 \end{itemize}
 Then $f$ has periodic orbits with homoclinic points.
 \end{prop}

We are now ready to show  theorem \ref{T_maim-e}.
\subsection{Proof of theorem \ref{T_maim-e}.}
 Let $\theta_t= \phi_t^\Omega(\theta) $ be a non-hyperbolic closed
orbit of minimal period $ T >0 $, contained in $ T^cM $. Let $ P=
P(\theta,\Sigma,\Omega) $ be the Poincar\'e map for a local
transversal section $ \Sigma \subset T^cM$ that contains  the
point $ \theta$. Since $ \theta_t $ is non-hyperbolic, the
eigenvalues of $ d_\theta P $ are of the form $ e^{\pm 2\pi i
\alpha} $, with $ \alpha \in [0,1)$. Recall that the
 symplectic twisted  form $ \omega(\Omega) $ induces a
 symplectic structure  in $ \Sigma $ and $P : \Sigma\rightarrow \Sigma $
preserves this  structure. Therefore, via  Darboux coordinates, we
can suppose that $ P $ is an area preserving diffeomorphism
defined in a neighborhood of $ 0 \in \R^2$ and $ P(0) = 0 $.

 Define $ Q\subset  J^3_s(1) $
 as:
$$ Q = \left\{\  \sigma\circ f_{\alpha,\beta}\circ \sigma^{-1}\ ;\  \sigma\in  J^3_s(1)
,\  \beta\not= 0,\ \mbox{ and }
% \ \alpha \notin \left\{ 0,
%\frac{1}{2}, \frac{2}{3},  \frac{4}{3},\frac{3}{2}\right\}
n\alpha \notin \N \mbox{ for all } n \in \{ 1,2,3,4\} \right\},$$
where $ f_{\alpha,\beta}: \R^2 \rightarrow \R^2 $ is given by $
f_{\alpha,\beta}( r,\theta) = ( r, \theta + \alpha + \beta r^2) +
\mathcal O (r^4) $, in polar coordinates.

 By the Birkhoff's normal form (theorem \ref{FNB}), the subset $ Q \subset J^3_s(1) $
is open and invariant. Since the orbit $\theta_t $ is
non-hyperbolic,  we have that $ j^3(P(\theta,\Sigma,\Omega)) \in
\overline Q$. Applying the theorem \ref{T-jet}, we obtain an exact
2-form $ d\overline\eta $ arbitrarily close to $ 0\in \Omega^2(M)
$ in the $C^r$-topology (with  $ r \geq 4 $)
 %y $ [\overline\Omega ] = [\Omega ] \in H^2(M,\R)$,
  such that
 $ \theta_t$ is a closed orbit of same period for the flow $
\phi_t^{\Omega+d\overline\eta} $\  and \
 $ j^k( P(\theta,\Sigma,\Omega+d\overline\eta) ) \ \in \ Q$.

Observe that $ \theta_t $ is  elliptic  for the perturbed flow
$\phi_t^{\Omega+d\overline\eta} $. Therefore, there is a
neighborhood $ \mathcal U \subset \Omega^2(M) $ of $ (\Omega+d
\overline \eta) $ such that, for all $ \overline \Omega \in
\mathcal U $, the flow $
\left.\phi_t^{\overline\Omega}\right|_{T^cM} $ has an
 elliptic closed orbit $ \overline \theta_t=\overline
\theta_t(\overline \Omega) $ close to $ \theta_t$  that we call
{\it analytic continuation of $ \theta_t$}. Since $ Q $ is open,
if the neighborhood  $ \mathcal U $ is taken small enough, we can
assume that $ j^3( P(\overline \theta, \Sigma, \overline\Omega))
\in Q $, for all $ \overline\Omega \in \mathcal U $.

\smallskip

On the other hand, by theorem \ref{T-KS}, there is a $
C^r$-residual subset $ \mathcal O(\Omega,c) \subset  \{
\overline\Omega \in \Omega^2(M) ; \ [\overline \Omega]=[\Omega]\}$
(with $r\geq 4$) such that, for each $ \overline \Omega \in
\mathcal O(\Omega,c) $, the corresponding magnetic flow $
\left.\phi_t^{\overline\Omega}\right|_{T^cM} $ satisfies the
property $ P_{K-S}$: all periodic orbits are elliptic or
hyperbolic and all heteroclinic orbits are transversal. Then, for
each closed orbit of $
\left.\phi_t^{\overline\Omega}\right|_{T^cM} $ its Poincaré map
satisfy  the conditions (i) and (ii) of proposition
\ref{p-calves}.

\smallskip

Since $ \m O(\Omega,c) $ is a residual subset of $ \{
\overline\Omega \in \Omega^2(M) ; \ [\overline \Omega]=[\Omega]\}$
we can $C^r$-approximate $d\overline \eta$ (with $r\geq 4$) by an
exact 2-form $ d\eta $ such that $ (\Omega +d\eta) \in \mathcal
O(\Omega,c) \cap \mathcal U $. Hence, if $ \overline \theta_t $ is
the analytic continuation of $ \theta_t $, then $f:= P(\overline
\theta, \Sigma, \Omega +d\eta ) $ satisfies  $ j^3(f) \in Q $
%and, via Darboux coordinates, $f$ is a diffeomorphism in a
%neighborhood of $ 0 \in \R^2 $ that preserves the area form $
%dx\wedge dy $
and the conditions (i) and (ii) of proposition \ref{p-calves}.

By definition of $ Q $, the map $ f $ is conjugated to a twist map
$  f_0 = h f h^{-1}$, in polar coordinates. In order to apply
proposition \ref{p-calves}, we need to  do a change of coordinates
which transforms  $ f_0 $ into a twist map  $ T : \R^+\times S^1
\rightarrow \R^+\times S^1 $, such that  the 1-form $ T^*(
Rd\theta) -Rd\theta $ is exact. Then the existence of a homoclinic
orbit implies the existence  of a non-trivial hyperbolic basic
set.

In fact, we consider  the following maps:
  $$
  \xymatrix{
  (x,y)\ar[r]&  (r,\theta)\ar[r]& (\frac{1}{2} r^2,\theta)=(R,\theta)
  \\
  \di  \ar[r]^{P\ \ }\ar[d]_{f}
  & \R^+\times S^1 \ar[r] \ar[d]_{f_0}& \R^+\times S^1\ar[d]^{T} \\
    \di \ar[r]& \R^+\times S^1 \ar[r]& \R^+\times S^1}
    $$
 where $ \di=\{ z\in \co ;\ |z|<1\ \}$,
 $P^{-1}(r,\theta)=(r\cos\theta, r\sin\theta)$.
 Let $G(x,y)=(\frac{1}{2} r^2,\theta)=(R,\theta)$.
 Then $\lambda:=G^*(R\ d\theta)=\frac{1}{2}(x\ dy-y\ dx)$.
Observe that $d\lambda=dx\wedge dy$ is the  area form $\di$.
 Since  $\di$ is contractible, we have that ${f_0}^*(\lambda)-\lambda $ is exact.
 Therefore $T^*(R\,d\theta)-R\,d\theta$ is exact. Since
 $R(r)=\frac {1}{2} \ r^2$
 strictly increasing  on $r>0$,  $T$
 is  a twist map if and only if  $f_{0}$ is a twist map.

$ \hfill{\Box}$

Let us  give two simple examples for  which we can apply the
theorem \ref{T_maim-e}.
\begin{example} \label{ej-S2}{\rm
Let $(M,g)$ be  a closed surface.
 We suppose  that the scalar  curvature
 satisfies $ \frac{1}{4} \leq K \leq 1 $. Let $ \Omega\equiv 0
\in \Omega^2(M)$. Then the $ \Omega$-magnetic geodesics are the
geodesics on $ M$  with  respect to the metric $ g$. In this case,
Thorbergsson proved in \cite{Thorb:p79}  the existence of a
non-hyperbolic closed geodesic.}
\end{example}

\begin{example}
\label{ej-moneda}{\rm Let $ B =\{ x\in \R^2, \|x\|^2<5 \}$, with
the  euclidian metric of $ \R^2$ and  the corresponding area  form
$ \Omega_0$. Let $ \eta $ be a 1-form in $ B $ such that $ d\eta =
-\Omega_0 $. We consider the exact magnetic field given by the
lagrangian
$$ L(x,v) = \frac{1}{2} \langle v,v \rangle - \eta_x(v). $$
 The Euler-Lagrange vector field of $ L:B \rightarrow \R $ can be seen as
local expression of a   magnetic field in a closed surface.

Since $ d\eta = - \Omega_0 $ the Euler-Lagrange vector field is
given by:
\begin{equation}
\label{ej-EL} \dot v = i\cdot v.
\end{equation}
We fix an  initial point $ p_0 = (-1,0 ) \in B $ and $ v_0=(1,0)
\in T_p^{\frac{1}{2}}B $. By the equation (\ref{ej-EL}), we have
that the circle   $ C:[0,2\pi]\rightarrow B $, given in polar
coordinates by $C(t)= (r(t),\theta(t))= (1,\pi-t)$ is a $
d\eta$-magnetic geodesics. Moreover, all  circles  obtained by
rotation of $ C $ fixing  the point $ p_0\in C $  are solutions of
(\ref{ej-EL}). Hence, if $ (p_0,v_0)\in \Sigma \subset
T_p^{\frac{1}{2}}B$ is a local transversal section, and $
P:\Sigma\rightarrow \Sigma $ is the corresponding Poincar\'e map,
then $ P( p_0,v) = (p_0,v ) ,\ \forall\  ( p_0,v ) \in \Sigma .$
Therefore, the orbit $ ( C(t),\dot C(t)) $ is degenerate, in
particular non-hyperbolic. }
\end{example}

\section{Franks' lemma for magnetic flows}
\label{Fraks-sec1} Let $ (M,g)$ be a closed riemannian surface and
let $ \nabla $ be the Levi-Civita connection. Then $ \nabla $
induces a connection $ K :TTM \rightarrow TM $, in the following
way: given $ \xi \in T_{\theta} TM $, let $ z : (-\epsilon ,
\epsilon ) \rightarrow TM $ be an adapted curve to $ \xi $ (this
is, $ z $ satisfies $ z(0) = \theta $ and $ \dot{z}(0)= \xi $),
then $ z(t) = ( \pi \circ z (t) , V(t) ) $, where V is a  vector
field  over $ \pi \circ z (t)$,  and we can define $ K_{\theta}
(\xi ) = \nabla_{\dot{ (\pi \circ z (t))}} V (0) .$ Let $
H(\theta), \ V(\theta) \subset T_\theta TM $ be the {\it vertical
and horizontal subspaces} defined as $$ V(\theta ) = Ker (
d_{\theta} \pi ) \ \mbox{ and }\ H(\theta) =
 Ker(K_{\theta})  $$
respectively. Then, we have the splitting:
$$ \Ttheta = H(\theta) \oplus V(\theta) \approx  T_{\pi(\theta)}M \times T_{\pi(\theta)}M. $$
Note that  the vertical subbundle does not depend on  the
riemannian metric and $ H(\theta ) $ and $ V(\theta) $ are
lagrangian subspaces of $( \Ttheta , \omega_0(\theta)) $.

\smallskip

Let us  recall the definition of the magnetic Jacobi fields. Let
$\Omega= f\Omega_0$ and $ \theta_t=( \gamma(t),\dot\gamma(t) ) $
be the orbit of a point $ \theta=(x,v) \in T^cM $, with $ c >0$.
Let $ z: (-\epsilon, \epsilon) \rightarrow TM $ be an adapted
curve of $ \xi \in \Ttheta $. We define  the {\it magnetic Jacobi
field } $ \J_\xi(t) $ as the vector field along  $ \gamma $ given
by:
$$  \J_\xi(t) = \left. \ps \right|_{s=0} \pi \circ \phi_t^\Omega(z(s)) .$$
Computing the horizontal and vertical components of the linearized
magnetic flow, we obtain:
\begin{equation}
\label{phi-j}
 d\phi_t^\Omega(\theta)(\xi) = \left( \J_\xi (t), \Dt \J_\xi (t) \right) \ \   ,\
 \forall \xi \in T_\theta TM.
\end{equation}
 %In fact, we have:
 %$$
 %\left\{ \begin{array}{cllr}
 %\J_\xi(t) &=& \left.\ps\right|_{s=0}( \pi \circ \phi_t(z(s)) =
 %  d_\theta(\pi \circ \phi_t)(\xi) = \Jxi ,\\
 %   \Dt \J_\xi (t)& =& \Dt \left.\ps\right|_{s=0} (  \pi \circ \phi_t(z(s)) =
 %    \left. \Ds \right|_{s=0} \pt ( \pi \circ \phi_t(z(s)) = \\
 %&=& \left. \Ds \right|_{s=0} \pt \gamma_s(t) = K_{\phi_t(\theta)}
 %( d \phi_t(\theta)(\xi))
 %\end{array} \right.$$
 Note that if $ \xi  \in T_\theta T^cM $  then $
 d\phi_t^\Omega(\theta)(\xi) \in  T_{\theta_t} T^cM $ for all $ t
 \in \R$. Therefore:
 \begin{equation}
 \label{derivada1}
 0\equiv dE  \left( \J_\xi(t) , \Dt \J_\xi(t)\right) = g\left( \dot
\gamma(t)\ ,\   \Dt \J_\xi(t) \right).
\end{equation}

Using that the curves $ t\mapsto \pi \circ \phi_t^\Omega(z(s)) $
are solutions of the equation (\ref{eq1}) and some basic
identities  of the Riemannian geometry, we have that a vector
field $ \J_\xi(t) $ along a $\Omega$-magnetic geodesic $ \gamma(t)
$ is a Jacobi field if, and only if, it satisfies the {\it
magnetic Jacobi equation} (for details see  \cite{P-P:c96}):
\begin{equation}
\label{eq-jacob} \Dtt \J_\xi + R(\dot\gamma, \J_\xi) \dot\gamma -
\left( \nabla_{\J_\xi} f \right)\ i \cdot \dot\gamma - f\ i \cdot
\Dt \J_\xi = 0,
\end{equation}
where $ R $ denotes the riemannian curvature tensor of $ (M,g) $.

 Let $ \mathcal N(t
)= \mathcal N(\theta_t) \subset T_\theta T^cM $ be the set
\begin{equation}
\label{bundle-N} \mathcal N (t) = \{ \xi \in T_{\theta_t} T^cM ;\
g_{\gamma(t)} ( d\pi(\xi) , \dot \gamma(t)) = 0 \} .
\end{equation}
 Since $ d_\theta \pi ( X^\Omega(\theta_t )) =
\dot\gamma(t) $, the subspace $ \mathcal N(t ) $ is transversal to
$ X^\Omega $ along   of $ \theta_t $. Hence
 $$ T_{\theta_t} T^cM = \mathcal N(t
 ) \oplus \langle X^\Omega(\theta_t)\rangle .$$
 % for all $ \theta \in T^cM $
  Therefore, the restriction of the  twisted  form $
 \omega_\Omega(\theta) $ to $ \mathcal N(\theta)$  is a non-degenerate 2-form.
 Note that $\mathcal N(\theta)$ does not
 depend on  $ \Omega$.
Let $ e_1(t) $ and $ e_2(t) $ be the vector fields along of $
\theta_t $ defined by:
\begin{equation}
\label{campos e-1-e-2} \left\{\begin{array}{cr}
e_1(t) = ( \ \ i\cdot \dot \gamma(t)\ ,\  0\  ) \mbox{ , }\\
e_2(t)= ( \ 0\ ,\  \ i\cdot\dot \gamma(t)\ )\ \ \
\end{array}\right. \in H(\theta_t) \oplus V
(\theta_t).
\end{equation}
%By the  definition of $ \mathcal N(t) $, we have
Note that $ e_1( t ), \
e_2(t) \in \mathcal N(t) $, %and
%$$ \omega_f({\theta_t}) ( e_1(t), e_2(t)) = g ( \ i \cdot
%\dot \gamma(t) , \ i \cdot \dot \gamma(t))= 2 c \not= 0 ,$$ for
for all $ t \in \R$.

 Let  $ \xi=\xi_1 e_1(0) + \xi_2e_2(0) \in T_\theta T^cM $. Using  the basis $
\{ \dot \gamma(t), i \cdot \dot \gamma(t) \}_{\gamma(t)} $  of $
T_{\gamma(t)} M $,  we can write $ J_\xi(t) = x(t) \dot\gamma(t) +
y(t) (\ i \cdot \dot\gamma(t))$ for some smooth functions  $ x,\ y
: \R \rightarrow \R $. A straightforward computation using
(\ref{phi-j}) and (\ref{derivada1}) shows that  $
x,y:\R\rightarrow \R  $ are solutions of
\begin{equation}
\label{eq-G-J} \left\{ \begin{array}{clcr}
\dot x &=& fy \hspace{4cm}\\
\ddot y &=& \{ -2cK_\gamma + g(\nabla f , \ i\cdot\dot \gamma)
-f^2 \} y \end{array} \right. ,
\end{equation}
with the initial conditions
 $ x(0) = 0 ,\  y( 0 ) = \xi_1 \ \mbox{ and } \ \dot y( 0 ) =\xi_2
 $,
and
\begin{eqnarray}
\label{flujo-m} d\phi_t^f(\theta)(\xi)& =& ( \J_\xi(t), \Dt
\J_\xi(t)) =
%\nonumber \\&=&
(  \ (x(t)\dot \gamma(t) + y (t)\ i\cdot \dot \gamma(t)) \  ,\
( \dot y (t) + x(t) f) \ i\cdot\dot \gamma(t) \  )=\nonumber\\
& = &  x(t) X^f(t) + y(t) e_1(t) + \dot y(t) e_2(t)
\end{eqnarray}

\smallskip

  Let $ \Sigma_t \in T^cM $ be a one parameter family of
  transversal sections such that $ T_{\theta_t} \Sigma_t  = \mathcal N(t) $.  Let \
   $ P_t(f): \Sigma_0 \rightarrow \Sigma_t $ be the    corresponding Poincar\'e
map, then
\begin{equation}
\label{dphi-M}
 d\phi_t^f(\theta)\left[\begin{array}{clcr}
0\\ \xi_1 \\ \xi_2 \end{array}\right] = \left[ \begin{array}{clcr}
 1 & \ * \\
0  & d_\theta P_t(f)
\end{array} \right]\left[\begin{array}{clcr}
0\\ \xi_1 \\ \xi_2 \end{array}\right].
\end{equation}

Using   ( \ref{flujo-m} ) and ( \ref{dphi-M} ), we obtain that $
d_\theta P_t(f): \mathcal N(0) \rightarrow \mathcal N(t) $ is the
fundamental matrix of
\begin{equation}
\label{eq-dP}
 \frac{d}{dt} \left[\begin{array}{cl} y(t)\\ \dot y(t) \end{array}\right] =
  \left[ \begin{array}{clcr}
0 & 1 \\
-\KM(f)(t) & 0
\end{array} \right]
\left[ \begin{array}{cl} y(t)\\ \dot y(t) \end{array}\right]
,\end{equation}
 where $ \KM(f)(t): = 2c K_{\gamma(t)} -
g_{\gamma(t)}(\nabla f, \ i \cdot \dot \gamma) + f^2(\gamma(t)) $
is called the  {\it magnetic curvature}.

\smallskip

 The following lemma is an easy consequence of the equation
(\ref{eq1}).

\begin{lem}\cite[lemma 2.1]{jan-generic:p2006}
\label{1-1} Given $ c > 0 $ and $ \Omega=f \Omega_0 \in
\Omega^2(M) $, let $ K= K(c,f) \in \R $ be  defined as
$K=\min\{\frac{1}{(\| f\|_{C^0}+1)^2}, \frac{i(M,g)}{2c}\}$, where
 $ i(M,g)$ denotes the injectivity radius. Then $ \pi \circ \phi^f ( \theta ) : [ 0 , K ) \rightarrow M $
 is injective, for all $ \theta \in T^cM$. In particular, any
 closed orbit of $ \phi_t^f $ in $T^cM$ has period less or equal to $K$.
\end{lem}

Let us now define our perturbation space. Let $ K=K(f,c) $ be
given by the lemma \ref{1-1}.\  For each $ \frac{K}{2} < T < K $
and $ W \subset M $ open with  $ W \cap \gamma((\frac{K}{2} ,T))
\not= 0 $, we define the set $ \mathcal{F}= \mathcal{F}(
W,\gamma,f,T)$ by:
$$ \mathcal{F}= \left\{ h \in C^\infty(M) ;
 \left.h\right|_{\gamma([0,T])} \equiv 0, \  \mbox{supp } ( h ) \subset W \mbox{ and } [h\Omega_0]=
  0 \in H^2(M,\R)\right\}  .$$
By definition, if $ h \in \mathcal F $,  the  curve $\gamma: [ 0,
T] \rightarrow M $ satisfies:
$$ \Dt \dot \gamma(t) = f(\gamma(t)) \ i\cdot \dot \gamma(t) =
 (f + h )(\gamma(t))\  \ i\cdot \dot \gamma(t)\ ,\ \forall t \in (0,T) .$$
Hence  $  [0,T]\mapsto( \gamma(t), \dot \gamma(t)) $ is an orbit
segment of the perturbed  magnetic flow $ \phi_t^{f+h}$, and we
can  define  the map:
\begin{eqnarray}
\label{S}
S_{T,\theta}: \mathcal{F } &\longrightarrow & Sp(1)\\
h &\longmapsto& d_\theta P_T(f+h) \nonumber.
\end{eqnarray}
\smallskip

We are now ready to state the analog of  the infinitesimal part of
the Franks' lemma \cite{Franks:p71}  for magnetic flows on
surface. In \cite{C-P:p2002}, G. Contreras and G. Paternain proved
a version of this lemma for geodesic flows on surfaces by
performing $ C^2 $-perturbations of the metric which correspond to
$ C^1 $-perturbations of the geodesic flow.

\bigskip

\begin{thm}
\label{franks-L} Let $ f_0  \in C^\infty(M) $, $ \frac{K}{2}< T
\leq K $ $( K =K(c,f_0) )$  and let $ \mathcal{U} $ be  a
neighborhood of $ f_0 $ in the $ C^1$-topology. Then  there is $
\delta = \delta(\mathcal U ,f_0,c)>0 $ such that the image of the
set  $ \mathcal U \cap \mathcal F(f_0,\theta,W,T) $ under the map
$ S_{T,\theta} $ contains a ball of radius $ \delta $ centered at
$ S_{T,\theta}(f_0)$. Moreover, if $ \gamma(t) $ is a closed
magnetic geodesic  of minimal period $ T_\theta $, then there is a
neighborhood $ U=U(\theta ,f_0,c,W,T_\theta)  \subset M $ of $
\gamma([T,T_\theta])$ such that the image of the set $ \mathcal U
\cap \{  f \in \mathcal F ; \mbox{ Supp}( f -f_0 )\subset W-U\} $
under the map $ S_{T,\theta} $ contains a ball of  radius $ \delta
$ centered at $ S_{T,\theta}(f_0)$.
 \end{thm}

\bigskip

We will  prove this theorem in the   subsection \ref{proba
T-franks-L}. For this, we  follow the strategy  in  the proof of
the Franks'  lemma for geodesic flows on surfaces  given by G.
Contreras and G. Paternain in \cite{C-P:p2002}.

We shall use theorem \ref{franks-L} on a finite number of segments
of a closed magnetic geodesic $ \gamma: [0,T_\theta] \rightarrow M
$ (where $ T_\theta $ denotes its minimal period) such that the
perturbations will be independent in each segment.
%Therefore, we
%going to produce a perturbation of the linearized Poincar\'e map.

Let $ \gamma(t) $ be a closed magnetic geodesic of minimal period
$ T_\theta$. We fix  $ t_0 \in ( \frac{K}{2}, K ] $ and $ n=
n(\theta)=n(T_\theta,t_0) \in \N $ such that $ T_\theta = n t_0 $.
For each $ 0 \leq i \leq n-1$, we define  $ \gamma_i = \gamma(
it_0 + t ) $. Then , given   a tubular neighborhood $ W_i \subset
M $ of $ \gamma_i((0,t_0))
   $,  we can define the map
$ S_{i,\theta} : \mathcal F ( f_0,\theta,W_i) \rightarrow Sp(1) $
by $S_{i,\theta}(f)= d_{(\phi^f_{it_0}(\theta))}  P( f
,\Sigma_{it_0},
 \Sigma_{(i+1)t_0})$,
for each $ i \in \{ 0, ...,n-1\}$.

     Let $ W_0 $ be tubular neighborhood of the segment  $ \gamma_0 $.
 Applying theorem \ref{franks-L} to  the map
$ S_{0, \theta}: \mathcal F ( f_0,\theta, c , W_0) \rightarrow
Sp(1) $ we obtain  $ \delta_0 >0 $ and  a neighborhood $ U_0 $ of
the curve $ \gamma_1
* ....* \gamma_{n-1}
$ as in the second  part of  theorem \ref{franks-L}. In the
following step, we need to choose  a tubular neighborhood $ W_1
\subset  U_0 $ and applying  again  theorem \ref{franks-L} to the
map  $ S_{1, \theta}: \mathcal F ( f_0,\theta, c , W_1)
\rightarrow Sp(1) $, we obtain  $ \delta_1 >0 $ and  a
neighborhood $ U_1 $ of $ \gamma_2
* ....* \gamma_{n-1} $. Proceeding in the same way,  we obtain $
\delta_i >0 $ and neighborhoods $ W_i $, for $ 0 \leq i \leq n-1
$. Hence, for    $ W = \cup_{i=0}^{n-1} W_i$, we can define
  \begin{eqnarray*}
\label{Sr}
 S_{\theta} : \mathcal F ( f_0,\theta,W) &\longrightarrow& \prod_{i=0}^{n-1}Sp(1)\\
 f \hspace{0.8cm} & \longmapsto&  \prod_{i=0}^{n-1}
 d_{(\phi^f_{it_0}(\theta))}  P( f ,\Sigma_{it_0}, \Sigma_{(i+1)t_0}) \nonumber
\end{eqnarray*}
   Applying  theorem \ref{franks-L} $n$ times we have proved:
\begin{cor}
\label{c-franks-L}Let $ \theta \in T^cM $ be  such that $
\phi_t^{f_0}(\theta )$ is a closed orbit. Let $ W, n(\theta) $ and
$ S_\theta $ be defined as above. Given a neighborhood $ \mathcal
U $ of $ f_0 $ in the $ C^1$ topology, there is  $ \delta=
\delta(f_0,\mathcal U, c)$ such that the image of the set $
\mathcal F(f_0, \theta, W ) \cap \mathcal U $ under the map $
S_\theta $ contains a product of balls of radius $ \delta $
centered at $( S_{0,\theta}(f_0),....,S_{(n-1),\theta}(f_0)) \in
\prod_{i=0}^{n-1}Sp(1) $.
\end{cor}

\subsection{Proof of theorem \ref{franks-L} } \label{proba
T-franks-L}

Let $ \psi: ( 0,T) \times (-\epsilon_0 , \epsilon_0) \rightarrow W
$ be  a  coordinate system  in $ W\subset M $  such that
 $ \psi( t , 0 ) = \gamma(t) $ and
  $ \left\{ \frac{\partial}{\partial t} ,
  \frac{\partial}{\partial x} \right\}_{\gamma(t)} = \{ \dot \gamma(t) ,
   i\cdot \dot \gamma(t)\}$,
for all $ t\in (0,T)$. Let    $ a: [-1,1] \rightarrow [-1,1] $ be
a smooth function  that satisfies:
\begin{eqnarray*}
\mbox{(i)} &  a(0)=0 & \mbox{(ii)} \  \mbox {  Supp}(a) \subset [-\frac{1}{2},\frac{1}{2}] \\
\mbox{(ii)} &    a^\prime(0)= 1 & \mbox{(iii)} \ \  \| a\|_{C^1} =1 \\
  \mbox{(iv)} &  \int_{-1}^1  a(x) dx = 0  .&
\end{eqnarray*}
We define  the subset $ \mathcal{H}$ of smooth  functions $ h :
(0,T) \times (-\epsilon_0 , \epsilon_0) \rightarrow \R $ such that
 %with Supp$( h ) \subset W $ and the function is the product function
$ h (t,x) = a_{\epsilon_0}(x) b (t),$ where $a_{\epsilon_0}(x) =
\epsilon_0\  a \left(\frac{x}{\epsilon_0}\right)$ and $ b \in
C^\infty(\R) $ with  Supp$(b) \subset ( 0,T) $.

\begin{lem}
\label{norma}
 Let $ h= a_{\epsilon_0}(x)b(t) \in \mathcal H $. Then:
\begin{eqnarray*}
& \mbox{{\rm (i)} }& [h\Omega_0] = 0, \\
 & \mbox{{\rm (ii)} }&  \overline f = f + h \in \mathcal F \ ,\ \ \forall h \ \in \mathcal H ,
  \mbox{{\rm and }} \ \forall f \in \mathcal F,\\\
& \mbox{{\rm (iii)} }& \| h \|_{C^1} \leq 2 \| b \|_{C^0} +
\epsilon_0 \| b \|_{C^1}.
\end{eqnarray*}
\end{lem}
{\bf Proof.} (i) Note that  if  $ \eta $ is a 1-form in $ M $ with
support in the neighborhood  $ W $ defined by $$ \eta|_W = \left(-
\int_{- 1}^x a_{\epsilon_0}(x)b(t)\
  ds\right) dt ,$$ then $ d\eta  = h \Omega_0 $.

(ii) It follows from (i) and the definition of the set $ \mathcal
H$.

(iii) If $ g:\R^2 \rightarrow \R $ is a function that has enough
differentiability, then:
$$ \| g \|_{C^1} \leq \sup_{x,t}| g(x,t) | +
\sup_{x,t}\left| \frac{\partial g(x,t)}{\partial x} \right| +
 \sup_{x,t}\left| \frac{\partial g(x,t)}{\partial t} \right|.$$
 Hence
 $$\| h \|_{C^1} \leq \sup_{x,t}| a_{\epsilon_0}(x)b(t) | +
  \sup_{x,t}| a_{\epsilon_0}^\prime (x) b(t) | + \sup_{x,t}| a_{\epsilon_0} (x)b^\prime(t) | .$$
 Since $ |a_{\epsilon_0}(x)| = \left| \epsilon_0
a(\frac{x}{\epsilon_0}) \right|\leq \epsilon_0 \| a \|_{C^0} \leq
\epsilon_0$, \ and \ $ | a_{\epsilon_0}^\prime(x) | = \left|
a^\prime(\frac{x}{\epsilon_0}) \right|\leq  \| a \|_{C^1} = 1$, \
for all $ x \in \R $, we have:
$$\| h \|_{C^1} \leq ( 1 + \epsilon_0) \| b \|_{C^0} + \epsilon_0 \| b \|_{C^1} \leq 2
 \| b \|_{C^0} + \epsilon_0 \| b \|_{C^1} .$$
$\hfill{\Box}$

We will  now fix some  constants and functions that  will be
useful in  the following lemma. By  changing  $ \mathcal U $ if
necessary, we can suppose that $ \| f \|_{c^0} \leq \|f_0 \|_{c^0}
+ 1 \ ,\ \ \forall f \in \mathcal U $. Therefore $ K(c,f_0) \leq K
( c,f) \ , \  \forall \ f \in \mathcal U $. Let us denote $
K(c,f_0)$ by $ k_0$  and set $ k_1= k_1(\mathcal U,c ) > 1 $, so
that, if  $  f \in \mathcal U $ and $ X(t)= X(f,\theta,t ) $ is a
fundamental matrix for  the equation ( \ref{eq-dP} ), then
\begin{equation}
\label{k_1} \| X(t) \| \leq k_1 \  \ \  \mbox{ and } \ \ \
\|X^{-1}(t)\|\leq k_1 \ \ ,\ \forall t \in [0,k_0]\ \mbox{ and } \
\forall \ \theta \in T^cM.
\end{equation}
Let $ 0 < \lambda << k_0/2 $ and \  $ k_2 = k_2(\mathcal U ,
\lambda, c) > 0 $ be
 such that:
\begin{equation}
\label{k_2} \max_{|t-k_0/2|\leq\lambda} \| X(t) - X(k_0/2)\| \leq
k_2 \ \ \mbox{ and } \ \ \max_{|t-k_0/2|\leq\lambda} \| X^{-1}(t)
- X^{-1}(k_0/2)\| \leq k_2\ ,
\end{equation}
for all $ f \in \mathcal U $ and $ \theta \in T^cM $.  If $
\lambda= \lambda( f_0, U ,c) $ is small enough, we have:
\begin{equation}
\label{k_2,k_1} 0 < k_2 < \frac{1}{16 \ k_1^3} < 1 < k_1.
\end{equation}
Let $ \delta_\lambda \ ,\ \Delta_\lambda : \R \rightarrow
[0,\infty) $ be  $ C^\infty $-approximations  of the
  Dirac delta  at the point $ \frac{k_0}{2} $, such that:
Supp$(\delta_\lambda ) \subset [ \frac{k_0}{2} - \lambda
,\frac{k_0}{2}) $, Supp$( \Delta_\lambda ) \subset (
\frac{k_0}{2}, \frac{k_0}{2}+ \lambda]$, \  $ \int \delta_\lambda
dt = \int \Delta_\lambda dt = 1 $ and   Supp$( \Delta_\lambda ) $
is an interval. Let $ k_3=k_3(\lambda) = k_3(f_0 ,\mathcal U , c)
$ be defined as:
\begin{equation}
\label{k_3}
 k_3 =  k_1^2 \left( \|\delta_\lambda \|_{C^0} + \|\delta_\lambda^\prime \|_{C^0} +
  \|\Delta_\lambda \|_{C^0} \| K_{mag}(f_0,c )\|_{C^0} + \frac{1}{2}\|\Delta_\lambda^{\prime\prime}
   \|_{C^0}\right) .
\end{equation}
Let $ 0 < \rho < 1/ (4k_1^2 k_3)$, by ( \ref{k_2,k_1} ) we have:
\begin{equation}
\label{rho} \frac{1}{k_1^2} - k_3\rho - 4k_1 k_2 > \frac{
1}{2k_1^2}.
\end{equation}
Finally, let $ \alpha: [0,k_0]\rightarrow [0,1] $ be a  $ C^\infty
$-approximation of the   characteristic function  of the set
$$ [0,T] - \left[ \gamma^{-1}\{\gamma((T,T_0))\} \cup \partial  \mbox{Supp}(\Delta_\lambda)
\right] $$ such that
\begin{equation}
\label{alpha} \int_0^T |\alpha(t) -1| dt \leq \rho.
\end{equation}

\begin{lem}
\label{cota} For an arbitrary  small parameter $ s $,  let $ h^s =
a_{\epsilon_0}(x) b^s(t)\ \in \mathcal H $ be such that $
b^{s=0}\equiv 0 $ and
\begin{equation} \label{del-b} \left.
\frac{\partial}{\partial s}\right|_{s=0} b^s(t) = \alpha(t)
\left\{ \delta_\lambda(t) a + \delta_\lambda^\prime(t) b - \left(
\Delta_\lambda(t) K_{mag}(f)(t) +
\frac{1}{2}\Delta_\lambda^{\prime \prime}(t)\right) c \ \right\}
\end{equation}
where $ a,b,c  \in \R $, and $ K_{mag}(f)(t) $ is the magnetic
curvature of $ f \in
 \mathcal U \cap \mathcal F $. Then
$$ \left\|\left. \frac{\partial}{\partial s} \right|_{s=0} S_{\theta, T} (f + h^s) \right\|
> \frac{1}{2k_1^3} \left\|\left[\begin{array}{cccc}
b&c\\
a& -b \end{array}\right]\right\|$$
\end{lem}
{\bf Proof.} Observe that:
\begin{eqnarray}
\label{form-k-mag} \KM(f+h^s)(t)& =& \KM(f)(t) - g_\gamma( \nabla
h^s ,\  i\cdot \dot \gamma) =\nonumber\\ &=& \KM(f)(t) -
\frac{\partial h^s(0,t)}{\partial x}
 %=\nonumber \\&=&
 =\KM(f)(t) -  b^s(t) .
\end{eqnarray}
% We denote  $ X_h(s,T) = d_\theta P_T(f+h^s) $.  From  (\ref{eq-dP}),
%the matrix    $ X_h(s,T) $ is a fundamental matrix in time $ t = T
%$ of the equation:
%\begin{equation}
%\label{dt-X}
% \pt X_h(s,t) = A^s(t)\cdot X_h(s,t)
%\end{equation}
%where $ A^s= \left[\begin{array}{clcr} 0& 1\\ -\KM(f+h^s) & 0
%\end{array}\right]= A^0+ \left[\begin{array}{clcr} 0& 0 \\ b^s(t)
%& 0 \end{array}\right]$.
%
%\smallskip
%
 We define $Z_h(T)= \left.\ds\right|_{s=0} d_\theta
P_T(f+h^s)$.% = \left.\ds\right|_{s=0}X_h(s,T) $. Since $X_h(s,t) $
%satisfies  the equation (\ref{dt-X}), we have:
%\begin{eqnarray*}
%\dt Z_h(t) & =&  \dt \left. \ps\right|_{s=0}X_h(s,t) = \left.\ps\right|_{s=0} \pt X_h(s,t) =\\
%& =& \left.\ps\right|_{s=0} \left(A^s(t)\cdot X_h(s,t)\right) =
% A^0 Z_h(t) +  \left[\begin{array}{clcr} 0& 0 \\
%\left.\ps\right|_{s=0} b^s(t) & 0 \end{array}\right] X^0(t),
%\end{eqnarray*}
%with $ X^0(t)= X_h(0,t)$.  Using the formula of   variation of
%parameters (c.f. \cite[pg. 48]{Hartm:b64}).
\  Since $ b^{s=0}\equiv 0 $, it follows from \cite[lemma
3.1]{jan-generic:p2006} that
%the equality    $ Z_h(0) = 0 $ that:
\begin{equation}
\label{dS-h}
 Z_h(T)=  X^0(T) \int_0^T \left(X^0(t)\right)^{-1}\left[\begin{array}{clcr} 0& 0 \\
  \left.\ps\right|_{s=0}b^s(t)& 0 \end{array}\right] X^0(t)\  dt.
\end{equation}
From  (\ref{del-b}) and integration by parts we obtain
\begin{eqnarray*}
 Z(T) &=& X(T) \left\{ \int_0^T
\alpha(t)\delta_\lambda(t) X^{-1}(t) \left[\begin{array}{cccc}
b&0\\a&-b \end{array}\right] X(t) \ dt + \right. \\
& + & \left.\int_0^T \alpha(t)\Delta_\lambda(t) X^{-1}(t)
\left[\begin{array}{cccc} 0&c\\0&0 \end{array}\right] X(t) \ dt
\right\}.
 \end{eqnarray*}
We will   write   $$ P(t)= \frac{1}{\alpha(t)}
\left(\left.\ds\right|_{s=0} b^s(t)\right) X^{-1}(t)
\left[\begin{array}{clcr} 0&0\\ 1 & 0
\end{array}\right]X(t)  ,$$  $$ Q_1(t) =
X^{-1}(t)\left[\begin{array}{clcr} b&0\\ a & -b
\end{array}\right]X(t) \ ,\ \ \ \  Q_2=
X^{-1}(t)\left[\begin{array}{clcr} 0&c\\ 0 & 0
\end{array}\right]X(t) $$ and $ Q_0 := Q_1 + Q_2 $.
 Then
 \begin{equation}
\label{P_Q}
 \int_0^T P(t)\ dt  = \int_0^T\delta_\lambda(t) Q_1(t) \ dt+ \int_0^T \Delta_\lambda(t)
  Q_2(t)\ dt .
\end{equation}
Using (\ref{k_1}), we obtain:
%\begin{equation}
%\label{Q_1}
$$\| \delta_\lambda(t) Q_1(t) \| \leq \|\delta_\lambda\|_{C^0} \| Q_1(t)\| \leq \|
\delta_\lambda\|_{C^0} \ k_1^2 \left\| \left[ \begin{array}{cccc}
b&0\\a&-b \end{array}\right] \right\|, \ \mbox{ and }
$$
$$\| \Delta_\lambda(t) Q_2(t) \| \leq \|\Delta_\lambda\|_{C^0} \| Q_2(t)\|
\leq \|\Delta_\lambda\|_{C^0} \ k_1^2 \left\| \left[
\begin{array}{cccc} 0&c\\0&0
 \end{array}\right]\right\| .$$
By (\ref{k_3}), we have:
\begin{eqnarray}
\label{P}
 \max_{t \in [0,T]} \| P(t) \|  &\leq&
   \left[ |a|\|\delta_\lambda \|_{C^0} + |b|\|\delta_\lambda^\prime \|_{C^0} +
  \left(\|\Delta_\lambda \|_{C^0} \| K_{mag}(f_0)\|_{C^0} +
  \frac{1}{2}\|\Delta_\lambda^{\prime\prime}
   \|_{C^0}\right) |c|\right]\  k_1^2 \nonumber \\
 &\leq &k_3  \left\|
  \left[ \begin{array}{cccc} b&c\\a&-b \end{array}\right]\right\|.
\end{eqnarray}
For each $ F: [0,T] \rightarrow \R^{2\times 2} $, we define:
$$ \mathcal O_\lambda\left( F,k_0/2\right) = \max_{|t-k_0/2|\leq \lambda}
 \left\| F(t) - F\left(k_0/2\right)\right\| \ \ . $$
Observe that if $ F,G : [0,T] \rightarrow \R^{2\times 2} $ and $ E
\in \R^{2\times 2} $ is a constant matrix  then:
\begin{eqnarray*}
\mathcal O_\lambda\left(F\ G, k_0/2\right)& \leq &
\max_{|t-k_0/2|\leq \lambda} \left\| F(t)G(t) - F\left(t \right)G
\left(k_0/2\right)\right\|+\\
&+&\max_{|t-k_0/2|\leq \lambda} \left\| F(t)G\left(k_0/2\right) -
F\left(k_0/2\right)G\left(k_0/2\right)\right\| \leq \\
&\leq&   \max_{|t-k_0/2|\leq \lambda} \left\| F(t)\right\|\mathcal
O_\lambda \left(G ,k_0/2\right) +\mathcal O_\lambda \left(F
,k_0/2\right)
% \max_{|t-k_0/2|\leq \lambda}\left\|G(t)\right\|
\left\|G\left(k_0/2\right)\right\| \\
\mathcal O_\lambda\left(E\ F, k_0/2\right) &=& \max_{|t-k_0/2|\leq
\lambda}
 \left\| E\ F(t) - E\ F\left(k_0/2\right)\right\| \leq \\
&\leq& \|E\|\  \mathcal O_\lambda\left(F, k_0/2\right).
\end{eqnarray*}
Write  $ A =  \left[ \begin{array}{cccc} b&c\\a&-b
\end{array}\right] $. By   (\ref{k_1}), (\ref{k_2}) and the
two inequalities above, we have that
\begin{eqnarray}
 \label{O_Q} &&\mathcal
O_\lambda\left(Q_0,k_0/2\right) = \mathcal
O_\lambda\left(X^{-1}(t) A X(t),k_0/2\right)\leq\hspace{3cm}
\nonumber \\
&\leq& \max_{|t-k_0/2|\leq
\lambda}\left\|X^{-1}(t)\right\|\mathcal
O_\lambda\left(AX(t),k_0/2\right) + \mathcal
O_\lambda\left(X^{-1}(t),k_0/2\right) \|A\|
\ \|X\left(k_0/2\right)\|\leq \nonumber \\
&\leq &\max_{|t-k_0/2|\leq \lambda}\left\|X^{-1}(t)\right\|
\|A\|\mathcal O_\lambda\left(X(t),k_0/2\right) + \mathcal
O_\lambda\left(X^{-1}(t),k_0/2\right) \|A\|
\ \|X\left(k_0/2\right)\|\leq  \nonumber\\
 &\leq& 2 k_1\ k_2 \| A\| .
\end{eqnarray}
By (\ref{P_Q}), we have
%\begin{eqnarray*}
$$ \left\|\int_0^T \alpha (t) P(t) \ dt - Q_0\left(k_0/2\right)\right\|
 =\left\|\int_0^T \alpha (t) P(t)dt -\int_0^T P(t) dt + \int_0^T
P(t) dt - Q_0\left(k_0/2\right)\right\| \leq
$$
\begin{eqnarray*}
& \leq& \left\|\int_0^T (\alpha(t)-1)P \ dt \right\| + \left\|
\int_0^T  \delta_\lambda (t) Q_1(t)-Q_1(k_0/2)+ \Delta_\lambda
(t)Q_2 - Q_2(k_0/2)\ dt
 \right\| \leq\\
& \leq & \int_0^T \| (\alpha(t)-1) P(t)\| \ dt + \mathcal
O_\lambda(Q_1,k_0/2) +
 \mathcal O_\lambda(Q_2,k_0/2) \leq\\
&\leq& \left(\max_{[0,T]}\|P(t) \|\right) \ \int_0^T |
(\alpha(t)-1)| dt + 2
 \mathcal O_\lambda(Q_0,k_0/2)).
 \end{eqnarray*}
Then   (\ref{P}), (\ref{O_Q}) and (\ref{alpha}) imply  that
\begin{equation} \label{cota1}
\left\|\int_0^T \alpha (t) P(t) \ dt -
Q_0\left(\frac{k_0}{2}\right)\right\|  \leq ( k_3\ \rho + 4\ k_1\
k_2 )  \|A\|.
\end{equation}
From the equality  $ A = X(k_0/2)\  Q_0(k_0/2)\  X^{-1} ( k_0/2 )
$ and (\ref{k_1}), we have
$$\|A\| = \| X(k_0/2)\
Q_0(k_0/2) \ X^{-1} ( k_0/2 ) \| \leq k_1^2 \|Q_0(k_0/2)\|\
\Longrightarrow$$
\begin{equation}
\label{|A|} \Longrightarrow \|  Q_0( k_0/2) \| \geq
\frac{\|A\|}{k_1^2}.
\end{equation}
Hence, using (\ref{cota1}), (\ref{|A|}) and (\ref{k_2,k_1})
\begin{eqnarray*}
  \left\| \int_0^T \alpha(t) P(t)\ dt \right\|& =& \left\| Q_0\left(k_0/2\right)+
  \int_0^T \alpha(t) P(t)\ dt  -  Q_0\left(k_0/2\right) \right\| \geq\\
& \geq & \| Q_0\left(k_0/2\right)\| - \left\|\int_0^T \alpha(t)
P(t)\ dt  -
 Q_0\left(k_0/2\right) \right\| \geq \\
&\geq & \left( \frac{1}{k_1^2} - k_3\ \rho - 4\ k_1\ k_2 \right)
\|A\| \geq \frac{1}{2k_1^2} \| A\|.
\end{eqnarray*}
Finally,
 the last inequality  and (\ref{k_1}) imply
$$ k_1\|Z(T) \| \geq \|X^{-1} Z(T)\| = \left\| \int_0^T \alpha(t)P(t) \ dt \right\| \geq
 \frac{1}{2k_1^2} \| A\|. $$
 Therefore $$ \| Z(T) \| \geq  \frac{1}{2k_1^3} \| A\|\ \ ,\ \forall f \in
\mathcal U \cap \mathcal F .$$ $ \hfill{\Box}$

We denote by $ \frak{sp} (1) $ the  Lie algebra  of the classical
Lie group $ Sp(1)= SL(2) $. For each matrix $ A=
\left[\begin{array}{cccc} b& c\\ a & -b
\end{array}\right] \ \in \frak{sp}(1) $,
let $ \beta_A: (0,T ) \rightarrow \R $ be defined as:
$$\beta_A(t) = \alpha(t)\left (\delta_\lambda(t) a +
 \delta_\lambda^\prime (t) b \right)
+ \left( K_{mag}(f_0)(t) + \frac{\Delta_\lambda^{\prime\prime}}{2
\Delta_\lambda(t)} \right)
  (e^{-\alpha(t) \Delta_\lambda(t) c }-1) .$$
We consider $ G: \frak{sp}(1)  \rightarrow C^\infty(M) $ such that
Supp$( G(A) ) \subset W $ and
\begin{equation}
\label{GA}
 \left. G(A)\right|_W =G(A)(t,x) = f_0(x,t) + a_{\epsilon_0}(x) \beta_A(t)
\end{equation}
in the tubular  coordinates
 $( t,x)$ in $ W $.
\begin{lem}
\label{G-A} For  $ \epsilon_0 $ small enough, there is $ \delta_1=
 \delta_1( \mathcal U, f_0, c )$, such that, if $ \|A \| < \delta_1 $ then
$ G(A) \ \in \ \mathcal U \cap \mathcal F $.
\end{lem}

{\bf Proof.} By lemma \ref{norma}, we have that $$ \| G(A) - f_0
\|_{C^1} = \| a_{\epsilon_0}(x) \beta_A(t) \|_{C^1} \leq 2
\|\beta_A\|_{C^0} + \epsilon_0\|\beta_A\|_{C^1}.$$ Let $ \epsilon
> 0 $ be such that  $ B_{C^1}( f_0, \epsilon )  \subset \mathcal U $,
where  $ B_{C^1}( f_0, \epsilon ) $ denotes the $ C^1 $-ball of
radius $ \epsilon $ centered at  $ f_0 $. Since $ \Delta_\lambda
> 0 $ for all
 $ t \ \in $ Supp$(\alpha)$,  there is  $ k_6 = k_6(\lambda,f_0)=
 k_6(f_0,  \mathcal U, c ) < \infty $ such that
 $$ k_6 = \max_{\|A\| \leq 1} \|\beta_A\|_{C^1} .$$
If  $ \epsilon_0 $ is small enough, we can suppose $ \epsilon_0 <
 \frac{\epsilon}{2k_6} $.

We consider  $ k_5 = k_5(\lambda, f_0 )= k_5(f_0, \mathcal U ,c )
$ given by:
$$ k_5 = \|\delta_\lambda\|_{C^0} + \| \delta_\lambda^\prime \|_{C^0} +
\left( \|K_{mag}(f_0)(t)\|_{C^0}  \|\Delta_\lambda(t)\|_{C^0} +
\frac{1}{2}
 \| \Delta^{\prime \prime}_\lambda \|_{C^0} \right)  e^{\|\Delta_\lambda\|_{C^0}}. $$
Observe that, if $ |c| \leq 1 $ then
$$ | e^{-\alpha \Delta_\lambda c} - 1 | \leq |c| \max_{|c| \leq 1 }
\left|\frac{\partial}{\partial c } \left( e^{-\alpha
\Delta_\lambda c} - 1 \right)\right| \leq |c| \Delta_\lambda
e^{\|\Delta_\lambda\|_{C^0}} $$
 and therefore
 %if $ |c| \leq 1 $
$$ \left|K_{mag}(f_0)(t) +  \frac{\Delta_\lambda^{\prime\prime}}{2 \Delta_\lambda(t)}
\right| |e^{-\alpha(t) \Delta_\lambda(t) c }-1|
 \leq |c|\left( |K_{mag}(f_0)(t) | \Delta_\lambda(t) +
  \frac{1}{2} | \Delta^{\prime \prime}_\lambda | \right)  e^{\|\Delta_\lambda\|_{C^0}}.$$
Hence, if $ \|A\| \leq 1 $ then
$$ \|\beta_A \|_{C^0} \leq k_5 \|A\| .$$

Choose $ \delta_1= \delta_1(f_0,\mathcal U , c )<1$ such that $
2k_5 \delta_1 < \frac{\epsilon}{2} $. Then
$$ \| G(A) - f_0 \| \leq 2 k_5 \|A\| + \epsilon_0 k_6 < \epsilon
,$$ for all  $\|A \| < \delta_1 $.

$\hfill{\Box}$

  A proof of  the following  general lemma can be seen in \cite{C-P:p2002}.

\begin{lem}
\label{gene}\cite[lemma 4.4]{C-P:p2002} Let $ N $ be a
n-dimensional smooth manifold and let $ F: \R^n \rightarrow N $ be
a sufficiently differentiable map such that
$$ \| d_xF \cdot v \| \geq a > 0 \, \ \ \forall\ (x,v) \in T\R^n \ \ \mbox{with}\
\| v \| \leq 1 \ \mbox{ and } \ \ \| x \| \leq r.$$ Then for all $
0 < b < a \ r $
$$ \{ p\in N ;\  d(p,F(0)) < b \} \ \subset
 \ F(\{ x\in \R^n ; \ \| x \| \leq \frac{b}{a}\}
.   $$

\end{lem}

 We consider a 3-parameter  family of maps
  $\{ G(A); \ A \in \frak{sp}(1) \} $ given in (\ref{GA}).  Observe that
$$ \frac{\partial\beta_A}{\partial a } = \alpha(t)\delta_\lambda(t) \ \ \mbox{ , }\ \ \
\frac{\partial\beta_A}{\partial b } =
\alpha(t)\delta_\lambda^\prime(t)  $$ and since $
\delta_\lambda(t) \ \Delta_\lambda (t) \equiv 0 $, we have
\begin{eqnarray*}
\frac{\partial\beta_A}{\partial c } &+&
\alpha(t)\Delta_\lambda(t)\beta_A(t) =
 - \alpha(t) \Delta\lambda(t) \left( K_{mag}(f_0)(t) +
 \frac{\Delta_\lambda^{\prime\prime}}{2\Delta_\lambda(t) } \right)
  e^{-\alpha(t) \Delta_\lambda(t) c} + \\
&+& \alpha(t) \Delta_\lambda(t) \left( K_{mag}(f_0)(t) +
\frac{\Delta_\lambda^{\prime\prime}}{2\Delta_\lambda(t) } \right)
( e^{-\alpha(t) \Delta_\lambda(t) c} - 1) =\\
&=&- \alpha(t) \left(  K_{mag}(f_0)(t) \Delta_\lambda(t) -
\frac{1}{2}\Delta_\lambda^{\prime\prime}(t) \right).
\end{eqnarray*}
Hence
$$ \frac{\partial\beta_A}{\partial c } = -\alpha(t)\left\{ \left[ K_{mag}(f_0)(t)+
\beta_A(t) \right] \Delta_\lambda(t) -
\frac{1}{2}\Delta_\lambda(t)^{\prime\prime} \right\}. $$ Then, by
(\ref{form-k-mag}), we have: $$ \frac{\partial\beta_A}{\partial c
} = -\alpha(t)\left\{ ( K_{mag}(G(A))(t) \Delta_\lambda(t) -
\frac{1}{2}\Delta_\lambda^{\prime\prime}(t)
  \right\}.$$
 Therefore, the directional derivative  of $ \frak{sp}(1) \ni A \mapsto \beta_A $
  satisfies (\ref{del-b}).

\smallskip

Consider $ F : \frak{sp}(1) \rightarrow Sp(1) $,  given by : $$
F(A) = S_{\theta,T} \circ G (A) .$$ It follows  from  lemmas
\ref{G-A} and \ref{norma} that there is  $ \delta_1 =
 \delta_1(\mathcal U,f_0,c) >0$ such that
$ G(A)\  \in \ \mathcal U \ \cap\  \mathcal F( f_0, \theta , W
,T_0) $, for all $ A\in \frak{sp}(1)$, with  $ \|A\| < \delta_1 $.
Since $ \alpha(t) =0 $ in a neighborhood of the  intersection
points of $ \gamma([0,T]) $ with  $ \gamma([T , T_0]) $,
 for  $ \epsilon_0 $  small enough,
there is  a neighborhood $ U $ of $ \gamma([T, T_0 ])$ such that
Supp$(G(A)) \cap U = \emptyset $, for all $ A \in \frak{sp}(1) $.

By  lemma \ref{cota}, we have that $$ \| d_B F \cdot A \| \geq
\frac{1}{k_1^3}>0\ ,\ \ \mbox{ for all } \ \| B\|< \delta_1 \
\mbox{ and } \ \|A\| \leq 1 $$ and applying    lemma \ref{gene} to
$ F $, with  $ r = \delta_1 $ and $ a= \frac{1}{2k_1^3}$, we
obtain
\begin{eqnarray*} B_{Sp(1)}(S_{\theta,T}(f_0), \delta_1/ 2k_1^3)&
\subset &  F \left( \{ A \in \frak{sp}(1) ; \ \|A\|
 < \delta_1 \} \right) =\\&=&
  S_{\theta,T } \left( \{G(A)\in \mathcal F; \ \|A\| < \delta_1 \} \right)
  \subset\\
& \subset & S_{\theta,T} \left( \mathcal U \cap \{ f \in \mathcal
F; \ \mbox{ Supp}(f-f_0 )
 \subset W-U \} \right)
\end{eqnarray*}
This inclusion proves the theorem.

%CORREGIDO ATE AQUI EM 22/03

\section{Magnetic flow with infinitely  many closed orbits in an energy level   }
\label{th-h} Let $( M,g)$ be a smooth closed and oriented
riemannian surface. For each $ c
> 0 $, let $ \mathcal R^1 (M,c ) $ be the  set of $ \Omega \in
\Omega^2(M)$ such that all closed orbits of $ \phi_t^\Omega $ in $
T^cM $ are hyperbolic endowed with the
 $ C^1$-topology. Given $ h \in H^2(M,\R)$, we define   $
\mathcal F^1(M,c,h)\subset  \mathcal R^1 (M,c )  $ as:
 $$ \mathcal F^1(M,c,h)= \{C^1\mbox{-interior of } \mathcal R^1
(M,c )\} \cap \{  [\Omega]=h \ \}.$$ Given $ \Omega \in \mathcal
F^1(M,c,h) $, let  $ Per(\Omega,c) \subset T^cM  $ be the union of
all  hyperbolic periodic orbits  of minimal period of  $
\left.\phi_t^\Omega\right|_{T^cM}$. By  definition,  $
\overline{Per(\Omega,c)} \subset T^cM $ is a compact and invariant
subset.

\smallskip

We recall  that a compact and invariant subset $ \Gamma \subset
T^cM$
   is a {\it hyperbolic set } if there exists a splitting
(continuous) of $ T_\Gamma (T^cM) = E^s \oplus E^u \oplus E^c ,$
where $ E^c = \langle X^\Omega\rangle $ and there are constants $
C > 0 $ and $ 0<\lambda < 1$, such that:
\begin{itemize}
\item[(a)]$d \phi_t^\Omega( E ^{s,u} ) = E^{s,u}$
\item[(b)] $| d_\theta\phi_t^\Omega( \xi ) | \leq C \lambda^t|\xi|\ , \ \
\forall \ t>0, \ \theta\in \Gamma ,\  \xi\in E^s, $
\item[(c)]$ |
d_\theta\phi_{-t}^\Omega( \xi ) | \leq C \lambda^t|\xi|\ , \ \
\forall \ t>0, \ \theta\in \Gamma ,\ \xi\in E^u.$
\end{itemize}
Using  Mañé's techniques  on  dominated sprinting \cite{Man:p82}
and a version of   Franks' lemma  for  magnetic flows (corollary
\ref{c-franks-L}), we will prove that:

\begin{thm}\label{t-per-hip}
  If $ \Omega \in \mathcal F^1(M,c,[\Omega])$, then $
\overline{Per(\Omega,c)} \subset T^cM $ is a hyperbolic set.
\end{thm}

A  hyperbolic set is called  {\it locally maximal}, if  there is
an open neighborhood $ U $  of $ \Gamma $, such that  $ \Gamma $
 is the maximal invariant  subset  of $ U $, i.e.,
 $ \Gamma = \bigcap_{t\in \R} \phi_t^\Omega(U).$
A {\it hyperbolic basic set } is a  maximal hyperbolic  set  with
a dense orbit and  we say  {\it non-trivial} when it is not a
single  closed orbit. It is well known  that a non-trivial
hyperbolic basic set has positive topological entropy (cf.
\cite{Bow:p72}).

Standard   arguments of dynamical systems \cite[\S 6]{Kat:b95}
imply that the set $ \overline{Per(\Omega,c)} $ is locally
maximal, for all $ \Omega \in \mathcal F^1(M,c,[\Omega])$. Hence,
it follows by Smale's Spectral Decomposition Theorem
\cite{smale:p67} ( see too \cite{Kat:b95}) that:
 \begin{cor}
 \label{c-per-hip}
Let $ c>0 $ and $ \Omega \in \mathcal F^1(M,c, [\Omega]) $. We
suppose  that the number of periodic orbits of minimal period for
the magnetic flow $
 \left.\phi_t ^\Omega\right|_{T^cM} $ is infinite. Then $
 \overline{Per(\Omega,c)} $ has a non-trivial basic set. In particular $
 \left.\phi_t ^\Omega\right|_{T^cM} $ has positive topological entropy.
\end{cor}

In the following two subsections we will recall some definitions
and statements that we will need to prove the theorem
\ref{t-per-hip}.

\subsection{Periodic sequences  of linear symplectic maps} \label{dom-sp}

 A linear map $
T: \R^n \rightarrow \R^n $ is {\it hyperbolic} if $ T $  has not
eigenvalue of norm equal to 1. The {\it stable and unstable
subspaces of $T$} are defined as
$$ E^s(T) = \left\{v \in \R^n;\ \lim_{n \rightarrow \infty }
T^n(v) =0 \right\}\ \mbox{ and } \  E^u(T) = \left\{v \in \R^n;\
\lim_{n \rightarrow \infty } T^{-n}(v) =0 \right\},$$
respectively.

 Let $ GL(n) $ be the group of linear  isomorphisms of $
\R^n $. We say that a sequence  $ \xi: \Z \rightarrow GL(n)$ is
{\it periodic} if there is $ n_0 \in \Z $ such that $ \xi_{i+n_0}
= \xi_i $, for all $ i \in \Z$. We say that a periodic sequence $
\xi $ is {\it hyperbolic} if  the linear map $ \prod_{i=0}^{n_0-1}
\xi_i $ is hyperbolic. In this case, we denote the stable and
unstable subspaces of $ \prod_{i=0}^{n_0-1} \xi_{j+i} $  by $
E_j^s(\xi) $ and $ E_j^u(\xi) $, respectively.

Given two families  of periodic sequences  $\xi^{(\alpha)} = \{
\xi^{(\alpha)}: \Z \rightarrow GL(n);\ \alpha \in \mathcal A \} $
and $ \eta^{(\alpha)}=\{ \eta^{(\alpha)}: \Z \rightarrow GL(n);\
\alpha \in \mathcal A \}$, we define:
$$ d(\xi^{(\alpha)},\eta^{(\alpha)}) =
\sup\{ \| \xi_n^{(\alpha)}- \eta_n^{(\alpha)}\| ;\
\alpha \in \mathcal A, \ n \in \Z \}.$$ We say that two periodic
families  in $ GL(n) $ are {\it periodically equivalent} if they
have the same index  set  $  \mathcal A $ and the minimal period
of $ \xi^{(\alpha)} $ and $\eta^{(\alpha)} $ coincide, for all $
\alpha \in \mathcal A $. We say that a family $ \xi^{(\alpha)} $
is a {\it periodic hyperbolic family} if  every  sequence  in $
\xi^{(\alpha)} $ is hyperbolic. Finally , we say that a periodic
hyperbolic sequence $ \xi^{(\alpha)} $ is {\it stably hyperbolic}
if there is $ \epsilon
> 0 $ such that  any family $ \eta^{(\alpha)} $ periodically equivalent
to $ \xi^{(\alpha)} $ satisfying  $ d
(\xi^{(\alpha)},\eta^{(\alpha)})< \epsilon $ is hyperbolic.

We are now ready to state the following result.
\begin{thm}{\rm (Mañé, \cite[lemma II.3]{Man:p82})}
\label{spl-dom} Let $\xi^{(\alpha)} = \{ \xi^{(\alpha)}: \Z
\rightarrow GL(n);\ \alpha \in \mathcal A \} $  be a stably
hyperbolic family of periodic sequences in $GL(n)$. Then there
exist constants $ m \in \N $ and $ 0<\lambda< 1 $, such that, for
all $ \alpha \in \mathcal A $ and $ j \in \Z $, we have:
$$ \left\|\left. \left(\prod_{i=0}^{m-1}
\xi_{i+j}^{(\alpha)}\right)\right|_{E_j^s(\xi^{(\alpha)})}
\right\| \cdot \left\| \left.\left(\prod_{i=0}^{m-1}
\xi_{i+j}^{(\alpha)}\right)^{-1}\right|_{E_{j+m}^u(\xi^{(\alpha)})}
\right\| \leq \lambda .$$
\end{thm}

 In \cite{C-P:p2002},
G. Contreras and G. Paternain proved  that if a family of periodic
hyperbolic sequences $ \xi^{(\alpha)} $ in $ Sp(1) $
 is stably hyperbolic among the  periodic sequences in $ Sp(1) $ and $ \sup_\alpha\|
\xi^{(\alpha)} \| < \infty $, then $ \xi^{(\alpha)}  $  is also
stably hyperbolic among the periodic sequences in $ GL(2)$.

\begin{cor}{\rm (\cite[corollary 5.2]{C-P:p2002})}
\label{c-spl-dom} Let $ \xi^{(\alpha)} = \{ \xi^{(\alpha)}: \Z
\rightarrow Sp(1);\ \alpha \in \mathcal A \}  $ be a family of
periodic hyperbolic sequences which is stably hyperbolic  in $
Sp(1) $, and $ \sup_\alpha\| \xi(\alpha)\| < \infty $. Then there
exist constants $ m \in \Z^+ $ and $ 0<\lambda< 1 $, such that,
for all $ \alpha \in \mathcal A $ and $ j \in \Z $, we have:
$$ \left\|\left. \left(\prod_{i=0}^{m-1}
\xi_{i+j}^{(\alpha)}\right)\right|_{E_j^s(\xi^{(\alpha)})}
\right\| \cdot \left\| \left.\left(\prod_{i=0}^{m-1}
\xi_{i+j}^{(\alpha)}\right)^{-1}\right|_{E_{j+m}^u(\xi^{(\alpha)})}
\right\| \leq \lambda .$$
\end{cor}

\begin{rem}
\label{rem-spl-dom} {\rm Let $ T_j^N =\prod_{i=0}^{N-1}
\xi_{i+j}^{(\alpha)} $. Using that $ \|AB\| \leq \|A\| \ \|B\| $,
for all $ A,B \in GL(n) $, we have that
$$ \left\|\left. T_j^{mN}\right|_{E_j^s(\xi^{(\alpha)})} \right\|
\cdot \left\|
\left.\left(T_j^{mN}\right)^{-1}\right|_{E_{j+mN}^u(\xi^{(\alpha)})}
\right\| \leq \lambda^N ,$$} for all $ N >1 $ and for all $ \alpha
\in \mathcal A$, $ j \in \Z$.
\end{rem}

\subsection{ Partially hyperbolic symplectic action}

 A
{\it symplectic vector bundle} $ \pi : E\rightarrow B $  is a
vector bundle such that the transition maps  preserve the
canonical symplectic structure  in the fibres $ \R^{2n} $. Let $
\Psi: \R \rightarrow Sp(E) $ be a continuous  action   that
preserves the fiber and satisfies $ \Psi_{s+t} = \Psi_s \circ
\Psi_t $. The action $ \Psi $ induces a flow $ \psi_t: B
\rightarrow B $ such that $\psi_t \circ \pi = \pi \circ \Psi_t$.

We say  that an action $ \Psi $ is {\it partially hyperbolic}, if
there is an invariant splitting $E=S \oplus U $, $ T >0 $ and
$0<\lambda <1 $, such that:
\begin{equation}
\label{dom-ac}
 \left\|\left.\Psi_T\right|_{S(b)} \right\|\cdot
\left\|\left.\Psi_{-T}\right|_{S(\psi_t(b))} \right\|\leq \lambda
\ ,\  \ \forall \ b \in B,
\end{equation}
and we say that $ \Psi $ is {\it hyperbolic} if there is a
splitting  $ E = E^s \oplus E^u $, $ C > 0$ and $ \lambda > 0 $,
such that:
\begin{eqnarray*}
&(1)& | \Psi_t( \xi ) | \leq C e^{-\lambda t}|\xi|\ ,\ \  \forall
\
t>0, \ \xi\in E^s ,\\
&(2)& | \Psi_{-t}( \xi ) | \leq C e^{-\lambda t}|\xi|\ , \ \
\forall \ t>0, \ \xi\in E^u.
\end{eqnarray*}
 The  domination  condition (\ref{dom-ac}),
implies that the decomposition $ E=  S \oplus U $ is continuous.
By definition, a hyperbolic action is  partially hyperbolic. The
converse is not true in general, but in the symplectic case we
have:
\begin{thm}{\rm(
\cite[corollary 1]{Contr:p2002})}
 \label{ph-h} A
symplectic partially hyperbolic action,  with $ \dim( S )= \dim (
U ) $ and compact base $ B $, is hyperbolic.
\end{thm}

\subsection{ Proof of theorem \ref{t-per-hip}}
\label{pruba T-per-hip}

We fix  $ c>0$. Let  $ U \subset T^cM $ be an open set and $
\mathcal R^1 (U) $ be the set  of  2-forms $ \Omega $ on $ M$ such
that all closed orbits of $ \phi_t^\Omega $  completely contained
in $ U $ are hyperbolic, endowed  with the  $ C^1$-topology. Given
$ h \in H^2(M,\R)$, we consider the subset
 $$ \mathcal F^1(U,h)= \{C^1\mbox{-interior of } \mathcal R^1(U )\}
  \cap \{ \Omega \in \Omega^2(M);\ [\Omega]=h \ \}.$$ Given $
\Omega \in \mathcal F^1(U) $, let $ Per(\Omega,U) $ be the union
of all periodic hyperbolic orbits  of  minimal period for the flow
$ \left.\phi_t^\Omega\right|_{T^cM}$  completely contained in $
U$.

The following proposition is a local version  that implies
theorem \ref{t-per-hip}.
\begin{prop}
\label{p-per-hip} Let $ \Omega \in \mathcal F^1(U,[\Omega])$. Then
$ \overline{Per(\Omega,U)} \subset T^cM $ is a  hyperbolic set.
\end{prop}

{\bf Proof.}
 For each $ \theta= (x,v) \in T^cM $, let $ \mathcal
N(\theta ) \subset T_\theta T^cM $ be  defined as:
$$\mathcal N (\theta) = \{ \xi \in T_{\theta} T^cM ;\
g_{x} ( d\pi(\xi) , v) = 0 \} .$$We recall  that
 $ T_\theta T^cM = \mathcal N(\theta
 ) \oplus \langle X^\Omega(\theta)\rangle $
  and  the restriction of the   twisted symplectic  form $
\omega_ \theta(\Omega) $ in $ \mathcal N(\theta)$  is
non-degenerate (section \ref{Fraks-sec1}). Let  $ K=K(\Omega,c )$
be given by lemma \ref{1-1}.\ Given $ \phi_t^\Omega(\theta) = (
\gamma(t),\dot \gamma(t) ) \in Per(\Omega,U) $, we denote its
minimal  period by $ T_\theta $. Let $ n= n(\theta,\Omega) \in \N
$ be such that $ T_\theta =  n \ t_0 $, for some $ t_0 \in (
\frac{K}{2}, K ] $. For each $ 0 \leq i \leq n-1$, we define  the
segment $ \gamma_i = \gamma( it_0 + t ) $.
 % By lemma \ref{1-1},  $\gamma_i:[0,t_0] \rightarrow U $ is injective.
  Let $ \Sigma_i \in
TU $ be the local transversal sections at $ \phi_{it_0}^\Omega
(\theta) $  such that $ T_{\phi_{it_0}^\Omega (\theta)} \Sigma_i =
\mathcal N(\phi_{it_0}^\Omega (\theta)) $. We denote by
$$ S_{i,\theta}= dP(\Omega, \Sigma_i,\Sigma_{(i+1)}):\mathcal
N(\phi_{it_0}^\Omega (\theta)) \rightarrow \mathcal
N(\phi_{(i+1)t_0}^\Omega (\theta))
$$ the linearized Poincar\'e maps.

 Since $  \Omega \in \mathcal F^1(U,[\Omega])$, there
is a $C^1$-neighborhood $ \mathcal U \subset \{ \overline\Omega
\in \Omega^2(M);\ [\overline \Omega ]= [\Omega]\} $ of $ \Omega $
such that each orbit of $ Per(\Omega,U) $ has a hyperbolic
analytic continuation, because otherwise we could produce a
non-hyperbolic orbit.
 Observe that if $
\overline \Omega \in \mathcal U $,  denoting by $\overline
\theta_t= \overline \theta_t(\overline \Omega) $ the analytic
continuation  of $ \theta_t =\phi_t^\Omega(\theta) $, then $
\overline \theta_t $ intersects the sections $\Sigma_i $, $ 0 \leq
i \leq n(\theta) $. Therefore,  we can  cut $\overline \theta_t$
into the same number of segments as $ \theta_t $. Hence the family
 \begin{equation}
\label{fam-per}
 \xi(\overline \Omega) = \left\{
S_{i,\overline \theta_t} (\overline \Omega)\ ;\ \ \
\phi_t^\Omega(\theta)\in Per(\Omega,U) \ \mbox{ and } \ 0\leq i
\leq n= n(\theta,\Omega) \right\},
\end{equation}
with $ \overline \Omega \in \mathcal U $  is a periodic equivalent
family.

The following lemma is a consequence of corollary
\ref{c-franks-L}.
\begin{lem}
\label{L-est-hip}
 If $\Omega \in \mathcal F^1( U,[\Omega])$, then $ \xi(\Omega ) $
 is  stably  hyperbolic.
\end{lem}
{\bf Proof.} Let $ \Omega \in \mathcal F^1( U,[\Omega])$  and let
 $ \mathcal U \subset \{ \overline\Omega \in \Omega^2(M);\
[\overline \Omega ]= [\Omega]\} $ a $C^1$-neighborhood of $ \Omega
$  such that, for all $ \overline \Omega \in \mathcal U $, the
family $ \xi(\overline \Omega ) $  is hyperbolic. Let $ \delta (
\Omega, \mathcal U , c )> 0 $ be given by corollary
\ref{c-franks-L}. We suppose  that the  family
$$ \xi(\Omega) = \{ S_i(\theta_t) = S_{i,\theta} : \mathcal N(\theta_{it_0})
\rightarrow \mathcal N(\theta_{(i+1)t_0});\ \theta_t \in Per(
\Omega,U) ,\  0 \leq i \leq n(\theta_t)\}
$$  is not stably hyperbolic. Then there exist a  closed orbit $
\theta_t \in Per( \Omega , U ) $ and a sequence of linear
symplectic maps
 $ \eta_i: \mathcal N(\theta_{it_0}) \rightarrow \mathcal
N(\theta_{(i+1)t_0}) $ such that:
$$ \| S_i(\theta_t) - \eta_i \| >\delta \ \ \mbox{ and } \ \
\prod_{i=0}^{n(\theta_t)} \eta_i \ \mbox{  is not hyperbolic}.$$
Observe that the  perturbation space  in the  Franks' lemma
(theorem \ref{franks-L}) preserves the selected orbit $ \theta_t
$. By corollary \ref{c-franks-L}, there is a 2-form $ \overline
\Omega \in \mathcal U $ such that $ \theta_t $  is a closed orbit
of $ \phi_t^{\overline\Omega} $
 too,   and $
 S_{i,\theta} (\overline \Omega) = \eta_i $. Since
 $$ d_{\theta_0}P(\overline \Omega, \Sigma_0,\Sigma_0) = \prod_{i=0}^{n(\theta_t)}
 S_{i,\theta} (\overline \Omega) =\prod_{i=0}^{n(\theta_t)} \eta_i
 , $$
 the orbit $ \theta_t $ is not hyperbolic for the magnetic flow
$ \phi_t^{\overline \Omega} $. This contradicts the choice of $
\mathcal U $.

 $\hfill{\Box}$

Applying corollary \ref{c-spl-dom} and remark \ref{rem-spl-dom},
if it is necessary, we have that there are $ 0<\lambda< 1 $ and  $
T
> 0 $ such that:
\begin{equation}
\label{dom-1} \left\| \left.d_{\theta} P_T \right|_{E^s(\theta)}
\right\| \cdot \left\|\left.d_{\theta_T} P_{-T}
\right|_{E^u(\theta_T)} \right\|  \leq \lambda,
\end{equation}
for all $ \theta_t = \phi_t^\Omega(\theta) \ \in Per(\Omega , U
)$, where $ P_T = P ( \Omega, \Sigma_0, \Sigma_T)$.

Let  $ \Gamma( U ) = \overline{Per(\Omega , U )} $. For each point
$ \theta \in \Gamma(U) $, we define:
$$ S(\theta):=\left\{ \xi \in \mathcal N (\theta)\  ; \ \begin{array}{clcr}
\exists \ \{\theta_n\}\ \subset \ Per(\Omega, U ),& \mbox{ with }
\lim_n \ \theta_n = \theta \mbox{ and }\\ \exists \ \xi_n  \in E^s
(\theta_n), & \mbox{ such that, } \lim_n\  \xi_n = \xi \end{array}
\right\}
$$
$$ U(\theta):=\left\{ \xi \in \mathcal N (\theta)\  ; \ \begin{array}{clcr}
\exists \ \{\theta_n\}\ \subset \ Per(\Omega, U ),& \mbox{ with }
\lim_n \ \theta_n = \theta \mbox{ and }\\ \exists \ \xi_n  \in E^u
(\theta_n), & \mbox{ such that, } \lim_n\  \xi_n = \xi \end{array}
\right\}.
$$
% $ \tau:\mathcal N \rightarrow \Gamma(U) $
Then the  uniform  domination condition (\ref{dom-1}) implies that
$ S \oplus U $ is a continuous splitting of $ \m N|_{\Gamma(U)}$
and  \begin{equation} \label{dom-2} \left\| \left.d_{\theta} P_T
\right|_{S(\theta)} \right\| \cdot \left\|\left.d_{\theta_T}
P_{-T} \right|_{U(\theta_T)} \right\| \leq \lambda,
\end{equation}
for all points $ \theta \ \in \Gamma(U) $.

Since $ dP_{s+t} = dP_s \circ dP_t$ and  (\ref{dom-2}),
 we have that $ dP:\mathcal N|_{\Gamma(U)}\rightarrow \mathcal
N|_{\Gamma(U)} $ is a  partially hyperbolic symplectic action. The
continuity  of the subbundles  $ S $ and $ U $ and their
definitions imply that $ S(\theta) = E^s(\theta) $ and $ U(\theta)
= E^u(\theta) $, for all $ \theta \  \in Per(\Omega,U) $. Hence $
\dim \ S = \dim U = 1$. By  theorem \ref{ph-h}, $ dP $ is a
hyperbolic symplectic action in $ \mathcal N|_{\Gamma(U)}$.

 Let $ \Lambda(\theta ) : \mathcal N(\theta ) \oplus
\langle X^\Omega(\theta) \rangle \longrightarrow \mathcal N(\theta
) $ be  the canonical projection. Observe that, by equality
(\ref{dphi-M}), we have:
$$d P  = \Lambda \circ \left.d\phi^\Omega \right|_{\mathcal N|_{\Gamma(U) }} .$$
Its follows from the graph transform  method of Hirsch-Pugh-Shub
in \cite{H-P-S:b77},
 that the hyperbolicity  of the action $ dP $ in $ \mathcal
N|_{\Gamma(U) } $ implies that $ \Gamma(U) $ is a  hyperbolic set
for the flow $ \phi_t^\Omega $.

$ \hfill{\Box}$

\subsection{ Proof of theorem \ref{T_maim-h}} Suppose $\Omega
\in \mathcal F^1(M,c,[\Omega]) $. By corollary \ref{c-per-hip},
the flow $ \left.\phi_t^{\Omega}\right|_{T^cM} $ has positive
topological entropy. Let $ \Omega \notin \mathcal
F^1(M,c,[\Omega]) $. Then there is $ d\overline \eta \in
\Omega^2(M)$ $C^1$-arbitrarily close to $ 0\in \Omega^2(M) $ such
that $ \left.\phi_t^{(\Omega +d\overline \eta)} \right|_{T^cM} $
has a non-hyperbolic closed orbit. Applying the theorem
\ref{c-per-hip}, we obtain an exact 2-form $ d \eta \in
\Omega^2(M)$,  $C^1$-arbitrarily close to $ 0\in \Omega^2(M) $,
such that  $\left.\phi_t^{(\Omega+d \overline \eta +
d\eta)}\right|_{T^cM} $ has positive topological entropy.

\section*{Acknowledgements}
This paper is a part of my PhD-thesis work held at CIMAT-M\'exico
under the guidance of the Profs. Renato Iturriaga and Gonzalo
Contreras. I am very grateful for their assistance. I thank
CAPES-Brazil and Fundep-UFMG for financial support. I would like
to thank the referees for the many comments and suggestions that
contributed for the improvement of this exposition.

\begin{appendix}
\label{Ap-lagr} \section{ Magnetic  Flows on the bidimensional
torus }\label{A-torus}

\subsection{Convex and superlinear autonomous lagrangians} \label{Ap-prel}

In this section, we will recall the main results for  convex and
superlinear autonomous lagrangian systems on closed manifolds such
as  the concepts introduced by Mather in \cite{Math:p91} and Mañé
in \cite{Man:c96}.

 Let $ M $ be a
closed smooth manifold  of arbitrary dimension. Let $ L : T
M\rightarrow \R $ be a smooth lagrangian that is  convex and
superlinear , i.e., for each fibre $ T_xM $, the restriction $
L(x,\cdot) $ has  positive definite  Hessian  and $
\lim_{\|v\|\rightarrow \infty} \frac{L(x,v)}{\|v\|} = \infty ,$
uniformly  in $ x \in M$. The {\it action} of $ L $ over an
absolutely continuous curve $ \gamma :[a,b]\rightarrow M $ is
defined as:
$$ A_L(\gamma) = \int_a^b L(\gamma(t),\dot \gamma(t).$$
The extremal  points  of the action are given by  solutions of the
{\it Euler-Lagrange equations} that in local coordinates can be
written  as:
$$ \dt \frac{\partial L}{\partial v } = \frac{\partial L}{\partial x
}.$$ Since the lagrangian $ L $ is convex, the Euler-Lagrange
equations define a complete vector field in $ T M $ that is
denoted by $ X_L $. We define  the {\it  Euler-Lagrange flow } $
\phi_t(L): T M \rightarrow T M $ as the flow associated to  $ X_L
$.

\smallskip

Let us  recall the main  concepts introduced by Mather in
\cite{Math:p91}, where  details and proofs (for periodic
lagrangians) can be found. Let $ \m M $ be the space of all
probability measure with compact support in the Borel
$\sigma$-algebra  of $ TM $ with the weak topology. %Basic
%results  on ergodic theory imply that $ \m M $ is a compact and
%convex metric space.
 We denote by $ \m M(L)\subset \m M $ the
subset of all invariant probability measures. We define the {\it
average  action } $ \m A_L: \m M (L) \rightarrow \R $  by $ \m
A_L(\mu) = \int_{TM} \ L \ d\mu. $ Given $ \mu \in \m M(L) $,
there is a unique homology class $ \rho(\mu) \in H_1(M,\R)$ such
that $ \langle \rho(\mu),[\omega]\rangle  = \int_{TM} \omega \
d\mu,$ for any closed 1-form $\omega $ on $ M$ . For each $ h \in
H_1(M,\R)$, we define $ \m M (h ) = \m M(L,h) = \{ \mu \in \m M(L)
\ ;\ \ \rho(\mu) = h \}. $ The {\it minimal action function} $
\beta:H_1(M,\R)\rightarrow \R $ is defined by $ \beta(h) = \min\{
\m A_L(\mu)\ ; \ \ \rho(\mu)=h\ \}.$ This function  is  convex and
superlinear. An invariant measure $ \mu $ that satisfies  $ \m
A_L(\mu) = \beta(\rho(\mu)) $ is called a {\it
$\rho(\mu)$-minimizing measure}. We define  the  $ \alpha $
 function as the  convex dual of $ \beta $, i.e. $ \alpha = \beta^*:
H^1(M,\R)\rightarrow \R $ is given by $ \alpha ([\omega]) =
\sup_{h \in H_1(M,\R)}\{<[\omega],h>- \beta(h) \}$. By convex
duality, we have that $ \alpha $  is also convex and superlinear,
and $\alpha^* = \beta$. Moreover, a measure $ \mu_0 $ is $
\rho(\mu_0)$-minimizing if and only if there is a closed 1-form  $
\omega_0 $, such that $ \m A_{L-\omega_0}(\mu_0) = \min_{\mu\in \m
M(L)}\{ \m A_{(L-\omega_0)}(\mu)\}$. Such a class $ [\omega_0]\in
H^1(M,\R) $ is  called a {\it subderivative } of
 $ \beta $  at the point $ \rho(\mu_0)$. For
each $ [\omega]\in H^1(M,\R)$, we consider  the subset:
\begin{equation}
 \label{conj-M-omega}
 \m M^\omega(L) = \{ \mu \in \m M (L)\ ; \ \ \m
\alpha([\omega])= -A_{L-\omega}(\mu )\ \}.
 \end{equation}
 {\it Mather's set}  $
\Lambda([\omega]) \subset TM $ is the compact and invariant subset
defined as:
$$ \Lambda([\omega]) = \overline{\bigcup_{\mu \in \m M^\omega(L)}\mbox{ Supp}( \mu )}
.$$ A important result of Mather in \cite[theorem 2]{Math:p91} is
that $ \Lambda([\omega])$  is a Lipschitz graph over a compact
subset of $ M$.

Using the properties of the $ \alpha $ and $ \beta $ functions, we
will prove:
\begin{lem}
\label{L-omega-en-c} We suppose  that $ H_1(M,\R)\not=0 $. Given $
c > c_0(L) $ and a nontrivial class  $ h_0 \in H_1(M,\R)$, there
is a closed 1-form  $ \omega_0 $ and $ \lambda_0 \in \R $, such
that $ [\omega_0] \in H^1(M,\R) $ is a subderivative of the  $
\beta $ function at the point $\lambda_0 h_0\in H_1(M,\R) $, with
$ \alpha([\omega_0]) = c $.
\end{lem}
{\bf Proof.} Since $ \beta $ is superlinear, we have:
\begin{equation}
\label{eq-sub-lambda}
 \lim_{\lambda \rightarrow  \infty}
\frac{\beta(\lambda h_0)}{|\lambda h_0|}= \infty.
\end{equation}
 Let
$ \partial \beta : H_1(M,\R ) \rightarrow H_1(M,\R )^* =H^1(M,\R )
$ be the multivalued function such that to each point $ h \in
H_1(M,\R ) $ associates  all subderivatives of  $ \beta $  in the
point $ h $.  Since $ \beta $ is finite,  $
\partial \beta(h) $ is  a  non empty  convex cone for all $ h
\in H_1(M,\R)$, and $
\partial \beta(h) $ is a unique vector if and only if $ \beta $ is differentiable
 in $ h $   ( cf \cite[Section 23]{Rock:b70}).
We define  the subset
$$ S(h_0)= \bigcup_{\lambda \in \R} \partial \beta (\lambda h_0 ). $$
By (\ref{eq-sub-lambda}) we have that the subset $ S(h_0) \subset
H^1(M,\R) $ is not bounded. Since $ \beta $ is continuous, by the
above properties of the multivalued function $ \partial \beta $,
we have that $ S(h_0) $ is a convex subset. Observe that if $
\omega \in
\partial \beta(0) $, then $\alpha([\omega]) = c_0(L) = \min
\{\alpha([\delta])\ ;\ \ \delta \in H^1(M,\R)\}  $ and, by
superlinearity of  $ \alpha $, the restriction $ \alpha|_{S(h_0)}
$ is not bounded. Therefore, by the  intermediate value theorem,
  for each $ c \in [c_0(L), \infty) $ there is $ \omega_0
\in \partial \beta (\lambda_0h_0)\subset S(h_0) $, for some $
\lambda_0 \in \R $, such that $ \alpha([\omega_0]) = c $.

$\hfill{\Box}$

Let us recall the definition of { \it Mañé's critical value}
 \cite{Man:c96,C-D-I:p97}. Given two points $ x,y \in M $ and $ T
> 0 $, we denote by $ \m C_T ( x,y ) $ the set of all
absolutely continuous curves $ \gamma:[0,T] \rightarrow M $, such
that $ \gamma(0) = x $ and $ \gamma(T) = y $. For each $ k \in \R
$, we define the { \it action potential  } $ \Phi_k: M \times M
\rightarrow \R $ by:
$$ \Phi_k(x,y) = \inf \{ A_{L+k}(\gamma)\ ;\ \ \gamma \in \cup_{T>0} \m C_T( x,y)\ \} $$
 The critical value  of a lagrangian $ L $ is the real number
 $ c(L) $ given by:
\begin{equation}
\label{c(l)}
 c(L) = \inf\{ k \in \R \ ; \ \  \Phi_k( x,x) > -\infty \
\mbox{ for any } x \in M \}.
\end{equation}

Let $ p: N\rightarrow M $ be a covering of $ M $ and $ \mathbb{L}:
TN \rightarrow \R $ be the lift  of  $ L $ to $ N $, i.e $ \mathbb
L = L \circ dp$. Then, for each $ k \in \R$, we can define the
action potential  $ \Phi_k: N\times N \rightarrow \R $ just as
above and we obtain the critical value $ c( \mathbb L ) $ for $
\mathbb L $. Among all  coverings of $ M $   the abelian covering,
which we shall  denote by $ \tilde M $, is distinguished. We
define the { \it strict critical value } $ c_0(L) $ as the
critical value (\ref{c(l)}) for the lift $\tilde L: \tilde
M\rightarrow \R $ of $ L $ to $ \tilde M $.

\smallskip

We recall that an absolutely continuous  curve $ \gamma:[a,b]
\rightarrow M $ is said to be  {\it semistatic} if $$
A_{L+c(L)}(\gamma|_{[t_0,t_1]}) =
\Phi_{c(L)}(\gamma(t_0),\gamma(t_1))\ \ \ \ \forall \ \  a <
t_0\leq t_1 < b .$$
 We say  that a semistatic curve $ \gamma|_{[a,b]} $  is
{\it static} when
 $ d_{c(L)}(\gamma(t_0),\gamma(t_1)) =0 $, for all $ a < t_0\leq t_1 < b $, where
 $ d_k( x,y ) =\Phi_k(x,y)  + \Phi_k (y,x)$ defines  a distance
function  for $ k > c(L) $ and a pseudo-distance  for $ k= c(L)$.
By definition of the action potential, a semistatic curve is a
solution of the Euler-Lagrange equations and has energy equal to $
c(L) $ (cf. \cite{Man:c96,C-D-I:p97}).
 The concepts of  semistatic  and static curves are
related   to the concepts of {\it  c-minimal trajectory  } and
{\it regular c-minimal trajectory   } respectively, that were
introduced  by Mather in \cite{Math:p93}.

\smallskip

Given $ \theta \in TM $ we will denote by $ \gamma_\theta:\R
\rightarrow M $ the unique solution of the Euler-Lagrange
equations with the  initial condition $
(\gamma_\theta(0),\dot\gamma_\theta(0)) = \theta$. We define the
sets:
$$ \Sigma = \Sigma(L) = \{ \theta \in TM \ ; \ \ \gamma_\theta:\R
\rightarrow M \ \mbox{ is semistatic }\},$$
$$ \widehat{\Sigma}=\widehat{\Sigma}(L)
 = \{ \theta \in TM \ ; \ \ \gamma_\theta:\R
\rightarrow M \ \mbox{ is static }\};$$ both  are compact and
invariant subsets of $ TM $. Let $ \pi:TM \rightarrow M $ the
canonical projection. Then $ \pi|_{\widehat{\Sigma}} :
\widehat{\Sigma} \rightarrow M $ is a bijective map with Lipschitz
inverse, the proof of this property can be seen  in \cite[theorem
VI]{C-D-I:p97} and is a extension of the Mather's graph theorem
\cite[theorem 2]{Math:p91}.

By the graph property, we can  define an  equivalence relation in
the set $ \widehat \Sigma $: two elements $ \theta_1 $ and $
\theta_2 \in \widehat \Sigma $ are equivalent when $ d_{c(L)}(
\pi(\theta_1),\pi(\theta_2))=0$. The equivalence relation breaks $
\widehat \Sigma $ into classes that are called {\it static classes
of L}. Let $ {\bf \Lambda} $ be the set of all static classes. We
define a partial order  $ \preceq $ in $ {\bf \Lambda} $ by: (i) $
\preceq $ is reflexive and transitive,
 (ii) if there is $ \theta \in \Sigma $, such that the  $
\alpha$-limit set  $ \alpha(\theta) \subset \Lambda_i$ and the  $
\omega$-limit set $ \omega(\theta)\subset \Lambda_j $, then $
\Lambda_i  \preceq \Lambda_j$. The following theorem was proved by
G. Contreras and G. Paternain in \cite[theorem A]{C-P:2002-3}.
\begin{thm}
\label{T-conect-stat}
 Suppose  that the number of static class is finite,
 then given $ \Lambda_i $ and $ \Lambda_j $ in $ {\bf \Lambda}$,
 we have that $ \Lambda_i \preceq \Lambda_j $.
 \end{thm}

Let $ \Lambda \subset TM $ be an invariant subset. Given $
\epsilon
> 0$ and $ T >0 $, we say that the points $ \theta_1, \theta_2
\in  \Lambda $ are  {\it $ (\epsilon,T)$-connected by chain} in
$\Lambda$, if there is a finite sequence $ \{ (\xi_i, t_i)
\}_{i=1}^n \subset \Lambda \times \R $, such that $ \xi_1=
\theta_1 $, \ $ \xi_n = \theta_2 $, \ $ T< t_i $ and $ d(
\phi_{t_i}(\xi_1), \xi_{i+1} ) < \epsilon$, for $ i = 1,..., n-1$.
We say that the subset $ \Lambda \subset TM $ is {\it chain
transitive}, if for all  $ \theta_1, \theta_2 \in \Lambda $ and
for all $ \epsilon >0 $ and $ T>0 $, the points $ \theta_1 $ and $
\theta_2 $ are $(\epsilon,T)$-connected by chain in $ \Lambda $.
When this condition holds for $ \theta_1=\theta_2 $, we say that $
\Lambda $ is {\it chain recurrent}. The proof of the following
theorem can be seen  in \cite[theorem V]{C-D-I:p97}.
\begin{thm}
\label{T-rec-estat}\
\begin{itemize}
\item[(a)] $ \Sigma $ is chain transitive.
\item[(b)] $ \widehat{\Sigma} $ is chain recurrent.
\item[(c)] the  $\alpha$ and $ \omega$-limit sets of a semistatic
orbit are contained in $ \widehat{\Sigma} $.
\end{itemize}
\end{thm}

We will now recall  the relations between  the Mather's theory and
the critical values of a lagrangian.   Mañé in
\cite{Man:c96,C-D-I:p97}, proved the notable equality $
c(L-\omega) = \alpha([\omega]), $ for all closed 1-form  $ \omega$
on $ M$.
 For the  abelian covering of $M$, in \cite{P-P:p97} G. Paternain and M.
Paternain proved :
$$ c_0(L)= \min_{[\omega]\in H^1(M,\R)}\{ c(L-\omega)\} =
\min_{[\omega]\in H^1(M,\R)}\{ \alpha([\omega])\}= -\beta(0).$$
 The following
theorem was also  proved  by Mañé in \cite{Man:c96}. The proof can
be see in \cite[theorem IV]{C-D-I:p97}.
\begin{thm}
\label{T-min-estat} Let $ \mu\in \m M(L)$ and $ \omega $ be a
 closed 1-form in $M$. Then $ \mu \in \m M^\omega(L) $ if only if
Supp$(\mu) \subset \widehat{\Sigma}(L-\omega)$.
\end{thm}

Finally, we recall that a class  $ h \in H_1(M,\R) $ is said to be
{\it rational } if there is $ \lambda>0$ such that $ \lambda h \in
i_*H_1(M,\Z) $, where $ i:H_1(M,\Z)\hookrightarrow H_1(M,\R) $
denotes the inclusion. The proof of the following proposition can
be found   in \cite{Mass:p97}. See also \cite[apendix
A]{C-M-p:p2004}.
\begin{prop}
\label{p-h-rational} Let $ M $ be a closed and oriented  surface.
Let $ \mu \in \m M(L) $ be a measure $ \rho(\mu)$-minimizing such
that $ \rho(\mu)$ is rational. Then the support of $ \mu $ is a
union of closed orbits of $ L $.
\end{prop}

\subsection{Magnetic flows on $T^2$} \label{L-mag-T2}

Let $ g $ be a riemannian metric in the torus $ T^2=\R^2/ \Z^2  $.
Given an exact 2-form  $ d\eta $ in $ T^2 $, the magnetic flow
associated to $ d\eta $ is given by the Euler-Lagrange equations
of the lagrangian $ L_\eta ( x,v) = \frac{1}{2} g_x(v,v) -
\eta_x(v) .$ The main result of this appendix is:
\begin{prop}
\label{pert-t2-infty} Given $ c >c_0(L_\eta)$
 there is an exact 2-form $ d\overline \eta \in \Omega^2(T^2) $
 arbitrarily $ C^r$-close to $ d\eta $, with $ 1\leq
r \leq \infty $, such that
 $h_{top}(d\overline \eta,c) > 0 $.
\end{prop}
{\bf Proof.}
 By  theorem \ref{T-KS}, there is a  $ C^r$-residual
subset $ \m O(c) \subset \{ \Omega \in \Omega^2(T^2)\ ;\ [\Omega]
=0 \} $ such that for all $ \Omega \in \m O(c) $ the exact
magnetic flow $ \phi_t^\Omega|_{T^cM} $  satisfies:  all closed
orbit are hyperbolic or elliptic   and all heteroclinic points are
transversal.
 Since $  L \mapsto c(L) $ is continuous in the set of the
lagrangians endowed with the uniform topology  in compact subsets
of $ TT^2 $ (cf. remark in \cite[pg17]{C-P:2002-3}), we can choose
an exact 2-form $ d\xi $ arbitrarily $ C^r$-close to $ d\eta $,
with $ 1\leq r \leq \infty$, that satisfies  $ d\xi \in \m O(c) $
and $ c
> c_0(L_{\xi})$.

Let $i: H_1(T^2,\Z) \hookrightarrow H_1(T^2,\R)$ be the inclusion.
Recall that $ H_1(T^2,\Z)=\Z^2 $ and that $ H_1(T^2,\R)= \R^2$.
Then $ \{(0,1),(1,0)\}\subset H_1(T^2,\Z) $ is a base of $
H_1(T^2,\R)$. It is easy to see that if $ \alpha_0,\alpha_1 $ are
two  closed curves in $ T^2 $, with $ [\alpha_0] =(0,1) $ and  $
[\alpha_1] =(1,0)$, then $ \alpha_0 \cap \alpha_1 \not= \emptyset
$.

By applying  lemma \ref{L-omega-en-c}\  for $ h_0= (0,1) \in
H_1(T^2,\Z) $, we obtain a closed 1-form  $ \omega_0 $ on $ T^2$
and $ \lambda_0 \in \R $, such that $ c= \alpha([\omega_0]) = c(
L_{\xi} -\omega_0) $. Hence,  if $ \mu_0 \in \m M(L_\xi) $ is a
$(\lambda_0 h_0)$-minimizing measure, then:
$$ \m A_{L_{\xi} - \omega_0}(\mu_0)= \min_{\mu\in \m M(L_\xi)}
\{ \m A_{L_{\xi} - \omega_0}(\mu)\}.$$ Since $ \rho(\mu_0) =
\lambda_0 h_0 $ is rational, by proposition \ref{p-h-rational},
 the support of $\mu_0 $ is a union of periodic
orbits.

Let $ \Lambda([\omega_0]) \subset TM $ the Mather's set associated
to  $ [\omega_0]\in H^1(T^2,\R) $. By definition of $
\Lambda([\omega_0]) $, we have that Supp$(\mu_0)\subset
\Lambda([\omega_0])$. By the graph theorem $
\pi|_{\Lambda([\omega_0])} : \Lambda([\omega_0])\rightarrow T^2 $
is injective. Hence, if $ \gamma_1, \gamma_2:\R\rightarrow T^2 $
are two distinct  closed $d\xi$-magnetic geodesics    contained in
$ \pi(\Lambda([\omega_0]))$, then  $ \gamma_1$ and $ \gamma_2 $
must be simple closed curves  and $ [\gamma_1]= n[ \gamma_2]\in
H_1(T^2,\Z)$, because otherwise  $\gamma_ 1 \cap \gamma_2 \neq
\emptyset $.
 Let $t\mapsto(\gamma(t),\dot \gamma(t))$
be   a closed orbit contained in the support of the measure $
\mu_0$ and let  $ \mu_\gamma $ be the  probability measure
supported in $ (\gamma(t),\dot \gamma(t))$. By definition of $
\rho $, we have that
$$ \rho(\mu_\gamma) = \frac{[\gamma]}{T} ,$$ where $ T
>0 $ denotes the period of $
\gamma $. It follows  from  $ \rho(\mu_0) = \lambda_0 h_0$ and
linearity of $ \rho$, that $ [\gamma]= n_0 h_0$ for some $0\not=
n_0 \in \Z$.

  Recall that  if a measure $\mu \in \m M^{\omega_0}(L_\xi) $,
 then all  ergodic components of $ \mu $
 are contained in $ \m M^{\omega_0}(L_\xi) $ too (cf. \cite[pg. 78]{Man:b91}).
Since $ c= c(L_\xi-\omega_0) > c_0(L_\xi)=-\beta(0) $ and  $ \m
M^{\omega_0}(L_\xi) $ is a closed set, there is $ k > 0 $ such
that
\begin{equation}
\label{0-fura-lambda}
 |\rho(\mu)|> k >0\ ,\ \ \mbox{ for all  } \ \mu\in \m
M^{\omega_0}(L_\xi) . \end{equation} Therefore, the period of a
periodic orbit contained in Supp$(\mu_0)$ is bounded. Then
Supp$(\mu_0)$ is a finite number of periodic orbits  of $ L_\xi $
(because $ d\xi \in \m O(c) $).

Let $ \mu \not= \mu_0 $ be  an invariant measure contained in $ \m
M^{\omega_0}$.
%Then by  graph property we have that
%Supp$(\mu)\cap $Supp$(\mu_0)= \emptyset$.
 We will  show that $
\rho(\mu) \in \langle h_0 \rangle \subset H_1(T^2,\R)$. Observe
that, if $ \gamma \in \pi($ Supp$(\mu_0))$ then $ [\gamma]\not=0$
and the set $ C_\gamma = T^2 - \{\gamma\} $ define an open
cylinder. By the
 graph property  of  $ \Lambda([\omega_0])$,
 we have that  Supp$(\mu)\cap $Supp$(\mu_0)= \emptyset$. Hence
  $\pi($Supp$(\mu)) \subset  C_\gamma$.
Therefore $ \rho(\mu ) \in i_* ( H_1(C_\gamma,\R)) \subset
H_1(T^2,\R)$, where $ i : C_\gamma \hookrightarrow T^2 $ is the
inclusion.  By definition of $ C_\gamma $ we have that $ \rho(\mu)
\in \langle h_0 \rangle \subset H_1(T^2,\R)$. Therefore
$$ \m M^{\omega_0}(L_\xi) \subset \{ \mu \in \m M(L_\xi) \ ;\ \ \rho(\mu)
\in \langle h_0\rangle \ \} .$$

 Using again  the proposition \ref{p-h-rational} and the inequality
(\ref{0-fura-lambda}), we obtain that the set $ \Lambda([\omega_0]
) $ is a  union of a finite number of periodic orbits for the
lagrangian $ L_{\xi} $.

If there is  a  non-hyperbolic closed orbit in $
\Lambda([\omega_0] ) $, then this proposition  reduces to theorem
\ref{T_maim-e}. It remains to consider the case  where all
periodic orbits  are hyperbolic. Let$ \gamma_i:\R \rightarrow T^2$
(with $ i=1,...,n $) be  closed magnetic geodesics  such that
$$ \pi\left( \Lambda([\omega_0])\right) = \bigcup_{i=1}^n \
\gamma_i.$$  Since Supp$(\mu_0) \subset \Lambda([\omega_0]) $, we
have that $ [\gamma_i] =n_0 h_0 = (0,n_0) \in H_1(T^2,\Z) $, for
all $ i \in \{1,...,n\}$. Given $ \theta \in TM $ we denote by $
\gamma_\theta:\R \rightarrow M $ the unique solution of the
Euler-Lagrange equations of $ L_\xi$, with the initial condition $
(\gamma_\theta(0),\dot\gamma_\theta(0)) = \theta$. Let
$$ \widehat{\Sigma}(\omega_0)=\widehat{\Sigma}(L_\xi-\omega_0)
 = \{ \theta \in TM \ ; \ \ \gamma_\theta:\R \rightarrow M \
\mbox{ is static for } (L_\xi-\omega_0)\},$$ and let $ {\bf
\Lambda} $ be the set of all static classes for the lagrangian $
L_\xi-\omega $. Recall that $ \Lambda([\omega_0]) \subset
\widehat{\Sigma}(\omega_0)$, and since each static class is a
connected set  (proposition 3.4 in \cite{C-P:2002-3}), for each
$1\leq i \leq n $, the orbit $(\gamma_i(t),\dot \gamma_i(t))$ is
contained in a static class. On the other hand,  theorem
\ref{T-min-estat}\ implies  that each static class contains at
least one of the periodic orbits in the  set $
\Lambda([\omega_0])$. Hence the number of  static classes is
bounded by $ n $.

 Suppose   that $ \Lambda([\omega_0]) =
\widehat{\Sigma}(\omega_0)$. Then each closed orbit in $
\Lambda([\omega_0]) $ is a static class. Let $
\Lambda_1,....,\Lambda_n $  be the static classes for the
lagrangian $ L_\xi -\omega_0 $. Applying the theorem
\ref{T-conect-stat}, we obtain that $ \Lambda_i \preceq \Lambda_j
$, for all $ i,j\in \{ 1,...,n\}$. In particular $ \Lambda_1
\preceq \Lambda_1 $. Therefore,  there is a point $ \theta \in
\Sigma( \omega_0)= \Sigma(L_\xi -\omega_0) $ such that the  $
\alpha$-limit set $\alpha(\theta) \subset \Lambda_1 $ and the  $
\omega$-limit set $ \omega(\theta)\subset \Lambda_1$. Since $ (
\gamma_1(t),\dot \gamma_1(t)) = \Lambda_1 $ is a hyperbolic orbit
of $ \phi_t^{d\xi}|_{T^cM} $ and that $ d\xi \in \m O(c) $, we
have that $ \Lambda_1 $ has a transversal homoclinic orbit $
\phi_t^{d\xi}(\theta) $. Then $ h_{top}( d\xi,c) > 0 $, which
proves the theorem in this case.

 We will now consider  the case  $ \Lambda([\omega_0]) \not=
\widehat{\Sigma}(\omega_0)$. For each  $ \theta \in
\widehat{\Sigma}(\omega_0) \setminus \Lambda([\omega_0])$, by
 the graph property  of the static set $
\widehat{\Sigma}(\omega_0)$, the magnetic geodesic  $
\gamma_\theta  : \R \rightarrow T^2 $ has no self-intersection
points and $ \gamma_\theta \cap \pi(\Lambda([\omega_0])) =
\emptyset$. Moreover, by theorem \ref{T-min-estat}, we have that
the  $\alpha$-limit and $\omega$-limit sets  are contained in the
Mather set $ \Lambda([\omega_0])$. Since    a curve on $ T^2$,
that  accumulates in positive time  to more than one closed curve,
must have self-intersection points, we have that  $ \omega(\theta)
$ is  a single closed orbit. By the same arguments, we have that $
\alpha(\theta) $ is  a single closed orbit. Since $d\xi \in \m
O(c) $, the orbit $ \phi_t^{d\xi}(\theta) $ is a transversal
heteroclinic orbit. Certainly, if $ \Lambda([\omega_0])$ is a
unique closed orbit, then $ \phi_t^{d\xi}(\theta) $ is a
transversal homoclinic orbit,  which implies   $
h_{top}(d\xi,c)>0$. In the opposite case, i.e, $ n
> 1$,
by recurrence property ( theorem \ref{T-rec-estat} ), we have that
if  $ \theta  \in \widehat{\Sigma}(\omega_0) \setminus
\Lambda([\omega_0])$, then
 $ \theta $ is an $(\epsilon,T)$-chain connected
in $\widehat{\Sigma}(\omega_0)$, for all $ \epsilon > 0$ and $ T
>0 $, i.e,  there is a finite  sequence $ \{ (\zeta_i, t_i)
\}_{i=1}^k \subset \widehat{\Sigma}(\omega_0) \times \R $, such
that $ \zeta_1= \zeta_k = \theta $, \ $ T< t_i $ and $ d(
\phi_{t_i}^{d\xi}(\zeta_1), \zeta_{i+1} ) < \epsilon$, for $ i =
1,..., k-1$.  Since the closed magnetic geodesics in $ \pi(
\Lambda([\omega_0]) )$ are isolated  on the torus, we have that
for $ \epsilon $ small enough, the set $ \{\pi(\zeta_i)\}_{i=1}^k
\subset \pi( \widehat{\Sigma}(\omega_0))$ must  intersect  the
interior of each one of the cylinders obtained by cutting  the
torus along the two  curves $ \gamma_i , \gamma_j \in \pi(
\Lambda([\omega_0]) )$, with $ 1\leq i,j \leq n $. Therefore,
choosing an orientation on $ \pi( \Lambda([\omega_0]) )$ and
reordering  the indices, we obtain a cycle  of   transversal
heteroclinic orbits. This implies that $ h_{top}(d\xi,c) >0$.

$\hfill{\Box}$
\end{appendix}

%-----------------------------------------------------------------
% ****************************************************************
% ----------------------------------------------------------------

\bibliographystyle{amsplain}

%\bibliography{XBib}
\def\cprime{$'$} \def\cprime{$'$}
\providecommand{\bysame}{\leavevmode\hbox
to3em{\hrulefill}\thinspace}
\providecommand{\MR}{\relax\ifhmode\unskip\space\fi MR }
\providecommand{\MRhref}[2]{%
  \href{http://www.ams.org/mathscinet-getitem?mr=#1}{#2}
} \providecommand{\href}[2]{#2}

\end{document}